%% file: annexe.tex
\documentclass[a4paper]{article}

\input francais_bis

\usepackage{amsmath}
\begin{document}

\newtheorem{lemme}{Lemme}
\newtheorem{defi}{D{\'e}finition}
\newtheorem{prop}{Proposition}
\newtheorem{pte}{Propri{\'e}t{\'e}}
\newtheorem{thm}{Th{\'e}or{\`e}me}
\newtheorem{cor}{Corollaire}

\def\dv#1#2{\langle {#1},{#2}\rangle}
\def\an#1#2{\def\deux{#2} \ifx\deux\empty {\cal O}_{#1} 
\else {\cal O}_{#1,#2} \fi }
\def\tk#1#2{{#2}\otimes _{#1}}
\def\tens{\raisebox{.3mm}{\scriptsize$\otimes$}}

\large
\large
\begin{center} \LARGE
{\bf Calculs explicites dans une alg{\`e}bre de Lie semi-simple effectu{\'e}s
  avec GAP4.}\\
\vspace{1cm}
{\sc Anne Moreau }\\

\vspace{1cm}

\footnotesize
\begin{center}
{\bf Abstract}
\end{center}
\begin{quotation}
\vskip .5em
In \cite{indice}, we show the following result, conjectured by D. Panyushev
\cite{Panyushev}, for $\g$ a semisimple Lie algebra\,:
\begin{eqnarray}\label{princ}
{\rm ind}~\n(\g^{e}) = {\rm rk}~\g-\dim \z(\g^{e}),
\end{eqnarray}
where $\n(\g^{e})$ and $\z(\g^{e})$ are, respectively, the normaliser
and the centre of the centraliser $\g^{e}$ of a nilpotent element
$e$. This result is proved in \cite{indice} when $\g$ is a classical
simple Lie algebra and when $e$ satisfies a certain property $(P)$. We
present in this paper the computations, made using GAP4, which prove that
distinguished, non-regular, nilpotent orbits in $E_6$, $E_7$, $E_8$ and
$F_4$ satisfy the property $(P)$. This work completes the proof,
presented in \cite{indice}, of the equality (\ref{princ}). The
complete proof of this result was already presented in \cite{indice_arxiv}.
\end{quotation}
\end{center}

\large

\section*{Introduction}\label{intro} 
Dans \cite{indice}, on prouve le r{\'e}sultat suivant, conjectur{\'e} par
D. Panyushev en \cite{Panyushev}\,:
\begin{thm}\label{intro1} Soit $e$ un {\'e}l{\'e}ment nilpotent d'une alg{\`e}bre
  de Lie semi-simple complexe $\g$. Alors on a\,: 
$${\rm ind}~\n(\g^{e}) = {\rm rg}~\g-\dim \z(\g^{e}),$$
o{\`u} $\n(\g^{e})$ et $\z(\g^{e})$ sont respectivement le normalisateur et
le centre du centralisateur $\g^e$ de l'{\'e}l{\'e}ment $e$.
\end{thm}
Ce r{\'e}sultat est d'abord prouv{\'e} dans \cite{indice} lorsque $\g$ est une alg{\`e}bre de Lie
simple classique (partie 3) et lorsque l'{\'e}l{\'e}ment $e$ v{\'e}rifie une certaine
propri{\'e}t{\'e} $(P)$ (partie 4), dont on rappelle ici la d{\'e}finition\,:
\begin{defi}\label{prop_P} Soit $\{e,h,f\}$ un
  $\mathfrak{sl}_2$-triplet dans $\g$ contenant $e$. On note $\z_{{\rm max}}$ le sous-espace propre de la restriction de ${\rm ad
} h$  {\`a} $\z(\g^{e})$ relativement {\`a} sa plus grande valeur
propre. On dira que $e$  {\rm v{\'e}rifie la propri{\'e}t{\'e} $(P)$}
si, pour tout {\'e}l{\'e}ment non nul $v$ de $\z(\g^{e})$, le sous-espace $\z_{{\rm max}}$ est contenu dans le sous-espace $[[f,\g^{e}],v]$.
 
Il est clair que si $e$ v{\'e}rifie la propri{\'e}t{\'e} $(P)$, il en est de
m{\^e}me de tous les {\'e}l{\'e}ments de l'orbite de $e$ sous l'action du groupe
adjoint. On dira qu'une orbite nilpotente de $\g$ {\rm v{\'e}rifie la
  propri{\'e}t{\'e} $(P)$} si l'un de ses repr{\'e}sentants la v{\'e}rifie.
\end{defi} 
En outre, il suffit de prouver le th{\'e}or{\`e}me \ref{intro1} pour les
{\'e}l{\'e}ments nilpotents distingu{\'e}s non r{\'e}guliers de $\g$. Ceci r{\'e}sulte
essentiellement de la proposition 2.4 de \cite{indice}. Il reste {\`a}
prouver dans \cite{indice}, que toutes les orbites nilpotentes distingu{\'e}es 
non r{\'e}guli{\`e}res d'une alg{\`e}bre de Lie simple exceptionnelles v{\'e}rifient
la propri{\'e}t{\'e} $(P)$. Il s'av{\`e}re qu'il suffit de tester la surjectivit{\'e} d'un nombre fini de matrices,
d{\'e}pendant d'un param{\`e}tre, pour v{\'e}rifier la propri{\'e}t{\'e} $(P)$. On
pr{\'e}sente ici les calculs effectu{\'e}s {\`a} partir de GAP4 qui permmettent de
v{\'e}rifier ces conditions. Ce travail ach{\`e}ve la d{\'e}monstration de la
proposition 5.3 de \cite{indice}. La preuve compl{\`e}te de ce r{\'e}sultat
est d{\'e}j{\`a} pr{\'e}sent{\'e}e dans \cite{indice_arxiv} .\\

On rappelle dans la partie \ref{rappel} les r{\'e}sultats de \cite{indice} que l'on justifie par les calculs de
GAP4. Dans la partie \ref{calcul}, on pr{\'e}sente ces calculs explicites. 

\section{Rappels des r{\'e}sultats de \cite{indice} que l'on justifie avec GAP4}\label{rappel}

On suppose que $\g$ est isomorphe {\`a} l'une des cinq
alg{\`e}bres de Lie simples exceptionnelles $E_6$, $E_7$, $E_8$, $F_4$ ou
$G_2$ et on suppose que $e$ est un {\'e}l{\'e}ment nilpotent distingu{\'e} non r{\'e}gulier de
$\g$. Il s'agit de montrer que l'{\'e}l{\'e}ment $e$ v{\'e}rifie la propri{\'e}t{\'e} $(P)$. On note
$m_{1},\ldots, m_{r}$ les valeurs propres de la restriction de ${\rm
  ad} h$ au sous-espace $\g^{e}$. Les entiers
$m_{1},\ldots,m_{r}$ sont pairs et on a 
$$2=m_{1} < m_{2} < \cdots < m_{r} \cdot$$
On note $\g^{e}_{m_{l}}$ le sous-espace propre correspondant {\`a} la valeur
propre $m_{l}$, pour $l=1,\ldots,r$. Avec les notations de la
d{\'e}finition \ref{prop_P}, on a\,: 
$$\z_{{\rm max}}=\g^{e}_{m_{r}} \cdot$$ 
On choisit une base 
$${\mathcal B}=
e_{m_{1}}^{1},\ldots,e_{m_{1}}^{d_{1}},e_{m_{2}}^{1},\ldots
,e_{m_{2}}^{d_{2}}, \ldots,e_{m_{r}}^{1},\ldots,e_{m_{r}}^{d_{r}}$$
de $\g^{e}$ de vecteurs propres telle que
$e_{m_{l}}^{1},\ldots,e_{m_{l}}^{d_{l}}$ forme une base de
$\g^{e}_{m_{l}}$, pour $l=1,\ldots,r$, et telle qu'il existe une base de
$\z(\g^{e})$ form{\'e}e de vecteurs de ${\mathcal B}$. On peut supposer que
$e_{m_{1}}^{1}=e_{2}^{1}=e$. 

Soit $i_{1} < \cdots < i_{s}$ dans
$\{1,\ldots,r\}$ et
$k_{(1,1)},\ldots,k_{(1,\delta_{1})},\ldots,k_{(s,1)},\ldots,k_{(s,\delta_{s})}$
des indices tels que les {\'e}l{\'e}ments 
$$e_{m_{i_{1}}}^{k_{(1,1)}},\ldots
,e_{m_{i_{1}}}^{k_{(1,\delta_{1})}},\ldots,
e_{m_{i_{s}}}^{k_{(s,1)}},\ldots
,e_{m_{i_{s}}}^{k_{(s,\delta_{s})}}$$  
forment une base de $\z(\g^{e})$. En particulier, on a les
relations\,: $m_{i_1}=2$ et $\delta_1=1$, $m_{i_s}=m_r$, $\delta_s=d_r$ et
$k_{(s,l)}=l$, pour $l=1,\ldots,\delta_s$.

Pour $i,j$ et $k$ dans $\{1,\ldots,r\}$ et $t,p$ et $q$ dans
$\{1,\ldots,d_{i}\}$, $\{1,\ldots,d_{j}\}$ et $\{1,\ldots,d_{k}\}$
respectivement, on note $\lambda_{(m_{k},q),(m_{i},t),(m_{j},p)}$ 
la coordonn{\'e}e de l'{\'e}l{\'e}ment $[[f,e_{m_{k}}^{q}],e_{m_{i}}^{t}]$ en
$e_{m_{j}}^{p}$. Notons que si
$\lambda_{(m_{k},q),(m_{i},t),(m_{j},p)}$ est non nul, on a la
relation $m_{k}=m_{j}-m_{i}+2$. 

Soit $l$ dans $\{1,\ldots,s\}$ et $\underline{\alpha_{l}}$  un
$\delta_l$-uplet. On d{\'e}finit une  matrice
$M(l,\underline{\alpha_{l}})$ de la fa{\c c}on suivante\,:\\ 
1) Si $m_{r}-m_{i_{l}}+2$ est une valeur propre
de la restriction de ${\rm ad} h$ {\`a} $\g^{e}$, on note $k(l)$
l'{\'e}l{\'e}ment de $\{1,\ldots,r\}$ tel que
$m_{k(l)}=m_{r}-m_{i_{l}}+2$. Alors $M(l,\underline{\alpha_{l}})$ est
la matrice de taille $d_{r} \times d_{k(l_0)}$
dont les coefficients $(M(l,\underline{\alpha_{l}}))_{p,q}$ sont donn{\'e}s par\,:
$$(M(l,\underline{\alpha_{l}}))_{p,q}
=\sum_{t=1}^{\delta_{l}} \alpha_{l}^{t}
\lambda_{(m_{k(l)},q),(m_{i_{l}},t),(m_{r},p)},$$
pour $p$ dans $\{1,\ldots,d_{r}\}$ et $q$ dans
$\{1,\ldots,d_{k(l)}\}$. La matrice $M(l,\underline{\alpha_{l}})$ s'{\'e}crit aussi, de mani{\`e}re
plus agr{\'e}able, comme une somme de matrices\,:
\begin{eqnarray}\label{M(v)}
M(l,\underline{\alpha_l})=\sum_{t=1}^{\delta_{l}} \alpha_{l}^{t}
M(l,\underbrace{
(0,\ldots,0,1,0,\ldots,0)}
_{{\rm avec \ } 
1 {\rm \ en \ } t^{{\rm i\grave{e}me}} {\rm \ position}}) \cdot
\end{eqnarray}
2) Sinon,  $m_{r}-m_{i_{l_0}}+2$ n'est pas une valeur propre de la restriction
de ${\rm ad} h$ {\`a} $\g^{e}$, et on pose\,:
$$M(l_0,\underline{\alpha_{l_0}})=0_{d_r,1} \cdot$$

La proposition suivante est d{\'e}montr{\'e}e dans \cite{indice}, Proposition 5.2\,:
\begin{prop}\label{M(v)_surj} On suppose que la
  matrice $M(l,\underline{\alpha_{l}})$ est surjective pour tout $l$ de
  $\{1,\ldots,s\}$ et tout $\delta_{l}$-uplet $\underline{\alpha_{l}}$
  non nul. Alors l'{\'e}l{\'e}ment $e$ v{\'e}rifie la propri{\'e}t{\'e} $(P)$.  
\end{prop}
{\bf Remarque} D'apr{\`e}s \cite{indice}, Lemme 5.1, il suffit de
  v{\'e}rifier la surjectivit{\'e} des matrices $M(l,\underline{\alpha_{l}})$, pour $l \geq 2$ et
  $\underline{\alpha_{l}}$ non nul. Cela laisse $s-1$ matrices,
  d{\'e}pendant d'un param{\`e}tre $\underline{\alpha_{l}}$, {\`a} {\'e}tudier. De
  plus, lorsque $\delta_s=d_r=1$, la matrice
$M(s,\underline{\alpha_s})=\alpha_{s}^{1} M(s,(1))$ est toujours
  surjective, pour $\alpha_{s}^{1}$ non nul. Dans ce cas, on r{\'e}duit {\`a}
  $s-2$ le nombre de matrices {\`a} {\'e}tudier.\\

Dans \cite{Triplets}, on trouve une liste de $\mathfrak{sl}_{2}$-triplets
correspondant aux orbites nilpotentes. Le logiciel GAP4 permet en
outre d'effectuer des calculs dans les alg{\`e}bres de Lie. 
Il permet notamment de calculer
le centralisateur d'un {\'e}l{\'e}ment, le centre d'une sous-alg{\`e}bre,
etc. En v{\'e}rifiant pour chaque orbite distingu{\'e}e non r{\'e}guli{\`e}re de
$E_6$, $E_7$, $E_8$ et $F_4$ les hypoth{\`e}ses de la proposition
pr{\'e}c{\'e}dente, on prouve la proposition suivante (Proposition 5.3 de \cite{indice}), puisque le cas de
$G_2$, particuli{\`e}rement simple, est trait{\'e} {\`a} part dans \cite{indice}\,:
\begin{prop}\label{verif} On suppose que $\g$ est une alg{\`e}bre de Lie
  simple de type $E_6$, $E_7$,
  $E_8$, $F_4$ ou $G_2$. Alors les orbites nilpotentes distingu{\'e}es non
  r{\'e}guli{\`e}res de $\g$ v{\'e}rifient la propri{\'e}t{\'e} $(P)$.
\end{prop}
Cette proposition termine la preuve du
th{\'e}or{\`e}me \ref{intro1}. Dans la partie
suivante, on expose les calculs qui justifient cette proposition. 

\section{Pr{\'e}sentation des calculs}\label{calcul}

L'alg{\`e}bre de Lie $\g$ est de type $E_6$, $E_7$, $E_8$ ou $F_4$. 
Pour une orbite nilpotente distingu{\'e}e non
r{\'e}guli{\`e}re de $\g$ donn{\'e}e, on consid{\`e}re l'{\'e}l{\'e}ment $e$ du
$\mathfrak{sl}_{2}$-triplet fournit par \cite{Triplets} correspondant {\`a}
la caract{\'e}ristique de l'orbite. Gr{\^a}ce {\`a} GAP4, on exhibe 
une base ${\mathcal B}$ v{\'e}rifiant les conditions pr{\'e}c{\'e}dentes. On v{\'e}rifie dans un premier temps que pour tout $l$ de $\{1,\ldots,s\}$, il
existe un entier $k(l)$ dans $\{1,\ldots,r\}$ tel que
$m_{k(l)}=m_{r}-m_{i_l}+2$. On calcule ensuite, la matrice
$M(l,\underline{\alpha_{l}})$, pour
$\underline{\alpha_{l}}=(\alpha_{l}^{1},\ldots,\alpha_{l}^{\delta_{l}})$
un $\delta_{l}$-uplet non nul, et on v{\'e}rifie la surjectivit{\'e} de cette
matrice. D'apr{\`e}s la remarque qui suit la proposition \ref{M(v)_surj}, on peut suppeoser
$l \geq 2$ et lorsque $d_r=1$, on peut supposer de plus $l \leq
s-1$. Enfin, dans la plupart des cas, on s'aper{\c c}oit que la somme de la
relation (\ref{M(v)}) n'a qu'un
seul terme. L'{\'e}tude de la surjectivit{\'e} de
$M(l,\underline{\alpha_{l}})$ ne d{\'e}pend
alors d'aucun param{\`e}tre, ce qui facilite le travail. 

La d{\'e}marche g{\'e}n{\'e}rale est la suivante\,: on d{\'e}finit l'alg{\`e}bre
de Lie {\tt L} dans laquelle on veut travailler gr{\^a}ce {\`a} la commande {\tt
  SimpleLieAlgebra}, on d{\'e}finit un syst{\`e}me de racines ({\tt
  RootSystem}), un syst{\`e}me de racines positives correspondant ({\tt
  PositiveRoots}) puis des syst{\`e}mes de vecteurs \guillemotleft
positifs\guillemotright \ et  \guillemotleft n{\'e}gatifs\guillemotright \ 
associ{\'e}s ({\tt PositiveRootVectors} et {\tt
  NegatitiveRootVectors}). La commande {\tt CanonicalGenerators} donne
une base de la sous-alg{\`e}bre de Cartan. On peut d{\'e}sormais faire des calculs dans
l'alg{\`e}bre de Lie {\tt L}. Il s'agit ensuite d'{\'e}tudier les orbites nilpotentes distingu{\'e}es non
r{\'e}guli{\`e}res. Pour chacune d'entre elles, on d{\'e}finit un
$\mathfrak{sl}_{2}$-triplet {\tt $\{$e,h,f$\}$} gr{\^a}ce aux donn{\'e}es de
\cite{Triplets}. On calcule ensuite le centralisateur {\tt g} de l'{\'e}l{\'e}ment positif {\tt e} avec la
commande {\tt LieCentralizer} puis le centre {\tt z} du centralisateur avec
{\tt LieCentre}. Pour chaque orbite, on pr{\'e}cise la valeur du plus haut
poids $m_{r}$ et on donne le nombre de matrices {\`a} {\'e}tudier. Pour
chacune d'entre elles on donne les valeurs de $m_{i_{l}}$ et de 
$m_{k(l)}$ et on effectue les calculs nec{\'e}ssaires. Les calculs de la premi{\`e}re orbite de la premi{\`e}re alg{\`e}bre (il s'agit de
l'orbite sous-r{\'e}guli{\`e}re de $E_{6}$) sont d{\'e}taill{\'e}s; les autres le sont
un peu moins.

\subsection{Calculs pour $E_{6}$}
\noindent On commence par d{\'e}finir {\tt L} et les g{\'e}n{\'e}rateurs de {\tt L}\,:
\begin{verbatim}
> L:=SimpleLieAlgebra("E",6,Rationals);
<Lie algebra of dimension 78 over Rationals>
> R:=RootSystem(L);
<root system of rank 6>
> P:=PositiveRoots(R);;
> x:=PositiveRootVectors(R);
[ v.1, v.2, v.3, v.4, v.5, v.6, 
v.7, v.8, v.9, v.10, v.11, v.12, 
v.13, v.14, v.15, v.16, v.17, v.18, 
v.19, v.20, v.21, v.22, v.23, v.24, 
v.25, v.26, v.27, v.28, v.29, v.30, 
v.31, v.32, v.33, v.34, v.35, v.36 ]
> y:=NegativetiveRootVectors(R);
[ v.37, v.38, v.39, v.40, v.41, v.42, 
v.43, v.44, v.45, v.46, v.47, v.48, 
v.49, v.50, v.51, v.52, v.53, v.54, 
v.55, v.56, v.57, v.58, v.59, v.60, 
v.61, v.62, v.63, v.64, v.65, v.66, 
v.67, v.68, v.69, v.70, v.71, v.72 ]
> CanonicalGenerators(R)[3]
[ v.73, v.74, v.75, v.76, v.77, v.78 ]
\end{verbatim}
Dans $E_{6}$, il y a deux orbites nilpotentes distingu{\'e}es non
r{\'e}guli{\`e}res\,: \\
\begin{enumerate} 
\item Caract{\'e}ristique\,: 
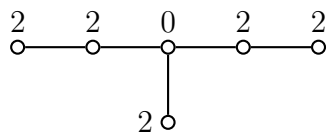
\begin{center}
\begin{pspicture}(-5,-1.2)(10,1)

\pscircle(-2,0){1mm}
\pscircle(-1,0){1mm}
\pscircle(0,0){1mm}
\pscircle(1,0){1mm}
\pscircle(2,0){1mm}
\pscircle(0,-1){1mm}

\psline(-1.9,0)(-1.1,0)
\psline(-0.1,0)(-0.9,0)
\psline(0.1,0)(0.9,0)
\psline(1.1,0)(1.9,0)
\psline(0,-0.1)(0,-0.9)

\rput[b](-2,0.2){$2$}
\rput[b](-1,0.2){$2$}
\rput[b](0,0.2){$0$}
\rput[b](1,0.2){$2$}
\rput[b](2,0.2){$2$}
\rput[r](-0.2,-1){$2$}
\end{pspicture}
\end{center}
On d{\'e}finit les {\'e}l{\'e}ments {\tt e} et {\tt f} du
$\mathfrak{sl}_{2}$-triplet correspondant dans les donn{\'e}es de \cite{Triplets}\,:
\begin{verbatim}
> e:=x[1]+x[2]+x[5]+x[6]+x[8]+x[9];
v.1+v.2+v.5+v.6+v.8+v.9
> f:=(12)*y[1]+(8)*y[2]+(-8)*y[3]+(22)*y[5]+(12)*y[6]+(8)*y[8]+
(22)*y[9]+(8)*y[10];;
\end{verbatim}
On v{\'e}rifie que le crochet {\tt e*f} est {\'e}gal {\`a} l'{\'e}l{\'e}ment neutre de
la caract{\'e}ristique et on pose {\tt h:=e*f}\,:
\begin{verbatim}
> e*f;
(12)*v.73+(16)*v.74+(22)*v.75+(30)*v.76+(22)*v.77+(12)*v.78
> h:=e*f;;
\end{verbatim}
On calcule le centralisateur {\tt g} de {\tt e} et on en donne une base
{\tt Bg}. On calcule ensuite le centre {\tt z} et on donne une base {\tt Bz} de {\tt z}\,:
\begin{verbatim}
> g:=LieCentralizer(L,Subspace(L,[e]));
<Lie algebra of dimension 8 over Rationals>
\end{verbatim}
Le centralisateur est de dimension 8.
\begin{verbatim}
> Bg:=BasisVectors(Basis(g));;
> z:=LieCentre(g);
<two-sided ideal in <Lie algebra of dimension 8 over Rationals>, 
(dimension 5)>
\end{verbatim}
Le centre est un id{\'e}al de dimension 5 dans {\tt g}.
\begin{verbatim}
> Bz:=BasisVectors(Basis(z));
[ v.1+v.2+v.5+v.6+v.8+v.9, v.23+(-1)*v.25+v.26, 
  v.27+(-1)*v.29+(-1)*v.30+(-1)*v.31, v.34+v.35, v.36 ]
\end{verbatim}
On calcule les \guillemotleft poids\guillemotright \  de {\tt z} en {\'e}valuant {\tt h*Bz[l]} pour
$l=1,\ldots,5$. On sait d{\'e}j{\`a} que {\tt h*Bz[1]=(2)*Bz[1]} car {\tt Bz[1]=e}.
\begin{verbatim}
> h*Bz[2];
(8)*v.23+(-8)*v.25+(8)*v.26
> h*Bz[3];
(10)*v.27+(-10)*v.29+(-10)*v.30+(-10)*v.31
> h*Bz[4];
(14)*v.34+(14)*v.35
> h*Bz[5];
(16)*v.36 
\end{verbatim}
On obtient que les poids sont $2,8,10,14,16$, d'o{\`u} $m_{r}=16$. Il y a
trois matrices {\`a} {\'e}tudier. 

\begin{enumerate}
\item $m_{i_{2}}=8$, $m_{k(2)}=10$. On cherche une base de $\g_{10}^{e}$
parmi les {\'e}l{\'e}ments {\tt Bg} en calculant {\tt h*Bg[i]}, pour $i=1,\ldots,8$ et on effectue le calcul correspondant\,:
\begin{verbatim} 
> h*Bg[5];
(10)*v.27+(-10)*v.29
> h*Bg[6];
(10)*v.30+(10)*v.31
\end{verbatim}
Le sous-espace $\g_{10}^{e}$ 
est engendr{\'e} par les vecteurs {\tt Bg[5]} et {\tt Bg[6]}.
\begin{verbatim}
> ((f*Bg[5])*Bz[2]);
(-20)*v.36
> ((f*Bg[6])*Bz[2]);
(20)*v.36
\end{verbatim}
La matrice {\`a} consid{\'e}rer est $M(2,(1))$; elle est
donn{\'e}e par $\left[
\begin{array}{cc}
-20 & 20
\end{array}
\right]$. C'est une matrice surjective.\\
\\
{\bf Remarque} Lorsque $d_{r}=1$, il suffit de trouver un {\'e}l{\'e}ment de $\g_{m_{k(l)}}^{e}$ qui donne un crochet non nul; dans la suite
on donnera seulement le calcul correspondant {\`a} cet {\'e}l{\'e}ment.

\item $m_{i_3}=10$, $m_{k(3)}=8$.
\begin{verbatim}
> ((f*Bz[2])*Bz[3]);
(-40)*v.36
\end{verbatim}

\item $m_{i_4}=14$, $m_{k(4)}=4$.
\begin{verbatim}
> h*Bg[2];
(4)*v.7+(2)*v.11+(2)*v.12+(2)*v.13+(-2)*v.14+(2)*v.15+(-4)*v.16

> ((f*Bg[2])*Bz[4]);
(14)*v.36
\end{verbatim}
\end{enumerate}
Ces trois calculs montrent que les hypoth{\`e}ses de la proposition \ref{M(v)_surj} sont v{\'e}rifi{\'e}es.\\
\\ 
{\bf Conclusion\,:} Cette orbite v{\'e}rifie la propri{\'e}t{\'e} $(P)$.\\

\item Caract{\'e}ristique\,: 
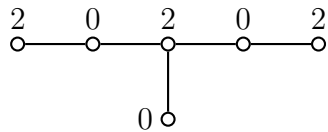
\begin{center}
\begin{pspicture}(-5,-1.2)(10,1)

\pscircle(-2,0){1mm}
\pscircle(-1,0){1mm}
\pscircle(0,0){1mm}
\pscircle(1,0){1mm}
\pscircle(2,0){1mm}
\pscircle(0,-1){1mm}

\psline(-1.9,0)(-1.1,0)
\psline(-0.1,0)(-0.9,0)
\psline(0.1,0)(0.9,0)
\psline(1.1,0)(1.9,0)
\psline(0,-0.1)(0,-0.9)

\rput[b](-2,0.2){$2$}
\rput[b](-1,0.2){$0$}
\rput[b](0,0.2){$2$}
\rput[b](1,0.2){$0$}
\rput[b](2,0.2){$2$}
\rput[r](-0.2,-1){$0$}
\end{pspicture}
\end{center}
D{\'e}finition du $\mathfrak{sl}_{2}$-triplet\,:
\begin{verbatim}
>e:=x[7]+x[8]+x[9]+x[10]=x[11]+x[19];
v.7+v.8+v.9+v.10+v.11+v.19
> f:=(8)*y[7]+(9)*y[8]+(5)*y[9]+(5)*y[10]+(8)*y[11]+y[19];
(8)*v.43+(9)*v.44+(5)*v.45+(5)*v.46+(8)*v.47+v.55
> e*f;                                                    
(8)*v.73+(10)*v.74+(14)*v.75+(20)*v.76+(14)*v.77+(8)*v.78
> h:=e*f;;
\end{verbatim}
Calcul de {\tt g}, {\tt Bg}, {\tt z} et {\tt Bz}\,:
\begin{verbatim}
> g:=LieCentralizer(L,Subspace(L,[e]));
<Lie algebra of dimension 12 over Rationals>
> z:=LieCentre(g);
<two-sided ideal in <Lie algebra of dimension 12 over Rationals>, 
  (dimension 4)>
> Bg:=BasisVectors(Basis(g));;
> Bz:=BasisVectors(Basis(z));
[ v.7+v.8+v.9+v.10+v.11+v.19, v.32+(-1)*v.33, v.35, v.36 ]
> h*Bz[2];
(8)*v.32+(-8)*v.33
> h*Bz[3];
(10)*v.35
> h*Bz[4];
(10)*v.36
\end{verbatim}  
Les poids de {\tt z} sont 2,8,10,10; d'o{\`u} $m_ {r}=10$. Il y a deux matrices {\`a} {\'e}tudier.

\begin{enumerate}
\item $m_{i_2}=8$,  $m_{k(2)}=4$
\begin{verbatim}
> h*Bg[4];
(4)*v.17+(-4)*v.18+(4)*v.20+(-4)*v.21
> h*Bg[5];
(4)*v.12+(4)*v.16+(-8)*v.22+(4)*v.24
> h*Bg[6];
(4)*v.22+(-4)*v.24+(-4)*v.25

> ((f*Bg[4])*Bz[2]);        
(-16)*v.35
\end{verbatim}

\item $m_{i_3}=10$,  $m_{k(3)}=2$
\begin{verbatim}
> h*Bg[2];
(2)*v.1+(2)*v.4+(2)*v.6+(2)*v.13+(2)*v.14+(-6)*v.15
> h*Bg[3];
(2)*v.19

> ((f*Bg[1])*Bz[3]); ((f*Bg[2])*Bz[3]); ((f*Bg[3])*Bz[3]);
(-10)*v.35
(10)*v.36
0*v.1
> ((f*Bg[1])*Bz[4]); ((f*Bg[2])*Bz[4]); ((f*Bg[3])*Bz[4]);
(-9)*v.36
(6)*v.35
(-1)*v.36
\end{verbatim}
La matrice {\`a} {\'e}tudier est
$$\alpha 
\left[
\begin{array}{ccc}
-10 & 0 & 0 \\
0 & 10 & 0 
\end{array}
\right] + \beta
\left[
\begin{array}{ccc}
0 & 6 & 0 \\
-9 & 0 & -1 
\end{array}
\right]=
\left[
\begin{array}{ccc}
-10\alpha & -\beta & 0 \\
-9 \beta & 10\alpha & -\beta 
\end{array}
\right]
 \cdot$$
On v{\'e}rifie que c'est une matrice de rang 2 pour tout
couple $(\alpha,\beta)$ non nul. 
\end{enumerate}
La proposition \ref{M(v)_surj} entraine que cette orbite v{\'e}rifie
$(P)$.\\
\\
{\bf Conclusion\,:} Cette orbite v{\'e}rifie la propri{\'e}t{\'e} $(P)$.\\
\end{enumerate}
{\bf Conclusion pour $E_6$\,:} Toutes les orbites nilpotentes distingu{\'e}es non
r{\'e}guli{\`e}res de $E_6$ v{\'e}rifient la propri{\'e}t{\'e} $(P)$.

\subsection{Calculs pour $E_{7}$}
\noindent D{\'e}finition de {\tt L}\,:
\begin{verbatim}
> L:=SimpleLieAlgebra("E",7,Rationals);
<Lie algebra of dimension 133 over Rationals>
> R:=RootSystem(L);
<root system of rank 7>
> P:=PositiveRoots(R);; 
> x:=PositiveRootVectors(R);
[ v.1, v.2, v.3, v.4, v.5, v.6, v.7, 
v.8, v.9, v.10, v.11, v.12, v.13, v.14, 
v.15, v.16, v.17, v.18, v.19, v.20, v.21, 
v.22, v.23, v.24, v.25, v.26, v.27, v.28, 
v.29, v.30, v.31, v.32, v.33, v.34, v.35, 
v.36, v.37, v.38, v.39, v.40, v.41, v.42, 
v.43, v.44, v.45, v.46, v.47, v.48, v.49, 
v.50, v.51, v.52, v.53, v.54, v.55, v.56, 
v.57, v.58, v.59, v.60, v.61, v.62, v.63 ]
> y:=NegativeRootVectors(R);
[ v.64, v.65, v.66, v.67, v.68, v.69, v.70, 
v.71, v.72, v.73, v.74, v.75, v.76, v.77, 
v.78, v.79, v.80, v.81, v.82, v.83, v.84, 
v.85, v.86, v.87, v.88, v.89, v.90, v.91, 
v.92, v.93, v.94, v.95, v.96, v.97, v.98, 
v.99, v.100, v.101, v.102, v.103, v.104, v.105, 
v.106, v.107, v.108, v.109, v.110, v.111, v.112, 
v.113, v.114, v.115, v.116, v.117, v.118, v.119, 
v.120, v.121, v.122, v.123, v.124, v.125, v.126 ]
> CanonicalGenerators(R)[3];
[ v.127, v.128, v.129, v.130, v.131, v.132, v.133 ]
\end{verbatim}
Dans $E_{7}$, il y a cinq orbites nilpotentes distingu{\'e}es non
r{\'e}guli{\`e}res\,:\\

\begin{enumerate}
\item Caract{\'e}ristique\,:
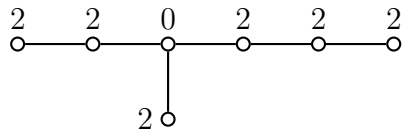
\begin{center}
\begin{pspicture}(-5,-1.2)(10,1)

\pscircle(-2,0){1mm}
\pscircle(-1,0){1mm}
\pscircle(0,0){1mm}
\pscircle(1,0){1mm}
\pscircle(2,0){1mm}
\pscircle(3,0){1mm}

\pscircle(0,-1){1mm}

\psline(-1.9,0)(-1.1,0)
\psline(-0.1,0)(-0.9,0)
\psline(0.1,0)(0.9,0)
\psline(1.1,0)(1.9,0)
\psline(2.1,0)(2.9,0)
\psline(0,-0.1)(0,-0.9)

\rput[b](-2,0.2){$2$}
\rput[b](-1,0.2){$2$}
\rput[b](0,0.2){$0$}
\rput[b](1,0.2){$2$}
\rput[b](2,0.2){$2$}
\rput[b](3,0.2){$2$}
\rput[r](-0.2,-1){$2$}
\end{pspicture}
\end{center} 
D{\'e}finition du $\mathfrak{sl}_{2}$-triplet\,:
\begin{verbatim}
> e:=x[1]+x[2]+x[3]+x[6]+x[7]+x[9]+x[11];
v.1+v.2+v.3+v.6+v.7+v.9+v.11
>f:=(26)*y[1]+(22)*y[2]+(50)*y[3]+(22)*y[5]+(40)*y[6]
(21)*y[7]+(15)*y[9]+(-15)*y[10]+(57)*y[11];
(26)*v.64+(22)*v.65+(50)*v.66+(22)*v.68+(40)*v.69
+(21)*v.70+(15)*v.72+(-15)*v.73+(57)*v.74
> e*f;
(26)*v.127+(37)*v.128+(50)*v.129+(72)*v.130
+(57)*v.131+(40)*v.132+(21)*v.133
> h:=e*f;;
\end{verbatim}
Calcul de {\tt g}, {\tt Bg}, {\tt z} et {\tt Bz}\,:
\begin{verbatim}
> g:=LieCentralizer(L,Subspace(L,[e]));
<Lie algebra of dimension 9 over Rationals>
> z:=LieCentre(g);
<two-sided ideal in <Lie algebra of dimension 9 over Rationals>, 
(dimension 6)>
> Bg:=BasisVectors(Basis(g));;
> Bz:=BasisVectors(Basis(z));
[ v.1+v.2+v.3+v.6+v.7+v.9+v.11, 
v.33+(2)*v.34+(-1)*v.36+(-1)*v.37+(-1)*v.38+v.40+(-3)*v.41, 
v.46+v.47+v.48+v.49+v.50, v.56+(-1)*v.57, v.60, v.63 ]
> h*Bz[2];
(10)*v.33+(20)*v.34+(-10)*v.36+(-10)*v.37
+(-10)*v.38+(10)*v.40+(-30)*v.41
> h*Bz[3];
(14)*v.46+(14)*v.47+(14)*v.48+(14)*v.49+(14)*v.50
> h*Bz[4];
(18)*v.56+(-18)*v.57
> h*Bz[5];
(22)*v.60
> h*Bz[6];
(26)*v.63
\end{verbatim}  
Les poids de {\tt z} sont 2,10,14,18,22,26; d'o{\`u} $m_{r}=26$. Il y a quatre matrices {\`a} {\'e}tudier.

\begin{enumerate}
\item $m_{i_2}=10$,  $m_{k(2)}=18$
\begin{verbatim}
> ((f*Bz[4])*Bz[2]);
(90)*v.63
\end{verbatim}

\item $m_{i_3}=14$,  $m_{k(3)}=14$
\begin{verbatim}
> ((f*Bz[3])*Bz[3]);
(-98)*v.63
\end{verbatim}

\item $m_{i_4}=18$,  $m_{k(4)}=10$
\begin{verbatim}
> h*Bg[3];
(10)*v.33+(20)*v.34+(-10)*v.36+(20)*v.37+(20)*v.38+(10)*v.40
> h*Bg[4];
(10)*v.37+(10)*v.38+(10)*v.41

> ((f*Bg[4])*Bz[4]);
(-30)*v.63
\end{verbatim}
\item $m_{i_5}=22$,  $m_{k(5)}=6$
\begin{verbatim}
> h*Bg[2];
(6)*v.19+(4)*v.20+(-2)*v.21+(2)*v.22
+(4)*v.23+(-2)*v.24+(-2)*v.25+(-2)*v.28

> ((f*Bg[2])*Bz[5]);
(-22)*v.63
\end{verbatim}
\end{enumerate}
{\bf Conclusion\,:} Cette orbite v{\'e}rifie la propri{\'e}t{\'e} $(P)$.\\

\item Caract{\'e}ristique\,:
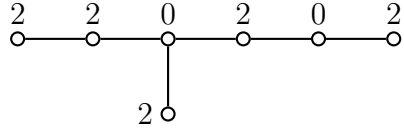
\begin{center}
\begin{pspicture}(-5,-1.2)(10,1)

\pscircle(-2,0){1mm}
\pscircle(-1,0){1mm}
\pscircle(0,0){1mm}
\pscircle(1,0){1mm}
\pscircle(2,0){1mm}
\pscircle(3,0){1mm}

\pscircle(0,-1){1mm}

\psline(-1.9,0)(-1.1,0)
\psline(-0.1,0)(-0.9,0)
\psline(0.1,0)(0.9,0)
\psline(1.1,0)(1.9,0)
\psline(2.1,0)(2.9,0)
\psline(0,-0.1)(0,-0.9)

\rput[b](-2,0.2){$2$}
\rput[b](-1,0.2){$2$}
\rput[b](0,0.2){$0$}
\rput[b](1,0.2){$2$}
\rput[b](2,0.2){$0$}
\rput[b](3,0.2){$2$}
\rput[r](-0.2,-1){$2$}
\end{pspicture}
\end{center} 
D{\'e}finition du $\mathfrak{sl}_{2}$-triplet\,:
\begin{verbatim}
> e:=x[1]+x[2]+x[3]+x[5]+x[7]+x[9]+x[18];                               
v.1+v.2+v.3+v.5+v.7+v.9+v.18
>f:=(22)*y[1]+(3)*y[2]+(42)*y[3]+(15)*y[5]+(17)*y[7]
+(28)*y[9]+(-28)*y[10]+(3)*y[12]+(3)*y[13]+(32)*y[18];
(22)*v.64+(3)*v.65+(42)*v.66+(15)*v.68+(17)*v.70
+(28)*v.72+(-28)*v.73+(3)*v.75+(3)*v.76+(32)*v.81
> e*f;
(22)*v.127+(31)*v.128+(42)*v.129+(60)*v.130
+(47)*v.131+(32)*v.132+(17)*v.133
> h:=e*f;;
\end{verbatim}
Calcul de {\tt g}, {\tt Bg}, {\tt z} et {\tt Bz}\,:
\begin{verbatim}
> g:=LieCentralizer(L,Subspace(L,[e]));
<Lie algebra of dimension 11 over Rationals>
> z:=LieCentre(g);                     
<two-sided ideal in <Lie algebra of dimension 11 over Rationals>, 
  (dimension 5)>
> Bg:=BasisVectors(Basis(g));;
> Bz:=BasisVectors(Basis(z));
[ v.1+v.2+v.3+v.5+v.7+v.9+v.18, v.39+(-1)*v.42+(-1)*v.43+(-1)*v.44
+(-2)*v.45+v.49, v.51+v.53+v.55+v.57, v.60, v.63 ]
> h*Bz[2];
(10)*v.39+(-10)*v.42+(-10)*v.43+(-10)*v.44+(-20)*v.45+(10)*v.49
> h*Bz[3];
(14)*v.51+(14)*v.53+(14)*v.55+(14)*v.57
> h*Bz[4];
(18)*v.60
> h*Bz[5];
(22)*v.63
\end{verbatim}  
Les poids de {\tt z} sont 2,10,14,18,22; d'o{\`u} $m_{r}=22$. Il y a trois matrices {\`a} {\'e}tudier.

\begin{enumerate}
\item $m_{i_2}=10$,  $m_{k(2)}=14$
\begin{verbatim}
> h*Bg[7];
(14)*v.53+(-14)*v.54
> h*Bg[8];
(14)*v.51+(14)*v.54+(14)*v.55+(14)*v.57

> ((f*Bg[7])*Bz[2]);
(-70)*v.63
\end{verbatim}

\item $m_{i_3}=14$,  $m_{k(3)}=10$
\begin{verbatim}
> h*Bg[5];
(10)*v.42+(10)*v.43+(10)*v.44+(20)*v.45
> h*Bg[6];
(10)*v.39+(10)*v.49

> ((f*Bg[5])*Bz[3]);
(70)*v.63
\end{verbatim}

\item $m_{i_4}=18$,  $m_{k(4)}=4$
\begin{verbatim}
> h*Bg[3];
(6)*v.20+(9)*v.21+(-3)*v.22+(-3)*v.27+(-9)*v.28
+(3)*v.29+(6)*v.30+(-3)*v.31+(-3)*v.35
> ((f*Bg[3])*Bz[4]);
(18)*v.63
\end{verbatim}
\end{enumerate}
{\bf Conclusion\,:} Cette orbite v{\'e}rifie la propri{\'e}t{\'e} $(P)$.\\

\item Caract{\'e}ristique\,:
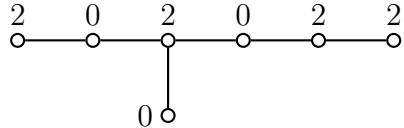
\begin{center}
\begin{pspicture}(-5,-1.2)(10,1)

\pscircle(-2,0){1mm}
\pscircle(-1,0){1mm}
\pscircle(0,0){1mm}
\pscircle(1,0){1mm}
\pscircle(2,0){1mm}
\pscircle(3,0){1mm}

\pscircle(0,-1){1mm}

\psline(-1.9,0)(-1.1,0)
\psline(-0.1,0)(-0.9,0)
\psline(0.1,0)(0.9,0)
\psline(1.1,0)(1.9,0)
\psline(2.1,0)(2.9,0)
\psline(0,-0.1)(0,-0.9)

\rput[b](-2,0.2){$2$}
\rput[b](-1,0.2){$0$}
\rput[b](0,0.2){$2$}
\rput[b](1,0.2){$0$}
\rput[b](2,0.2){$2$}
\rput[b](3,0.2){$2$}
\rput[r](-0.2,-1){$0$}
\end{pspicture}
\end{center} 
D{\'e}finition du $\mathfrak{sl}_{2}$-triplet\,:
\begin{verbatim}
> e:=x[7]+x[8]+x[9]+x[10]+x[11]+x[12]+x[22];                            
v.7+v.8+v.9+v.10+v.11+v.12+v.22
> f:=(15)*y[7]+(18)*y[8]+(24)*y[9]+(15)*y[10]+(10)*y[11]+
(28)*y[12]+y[22];
(15)*v.70+(18)*v.71+(24)*v.72+(15)*v.73+(10)*v.74+(28)*v.75+v.85
> e*f;
(18)*v.127+(25)*v.128+(34)*v.129+(50)*v.130
+(39)*v.131+(28)*v.132+(15)*v.133
> h:=e*f;;
\end{verbatim}
Calcul de {\tt g}, {\tt Bg}, {\tt z} et {\tt Bz}\,:
\begin{verbatim}
> g:=LieCentralizer(L,Subspace(L,[e]));
<Lie algebra of dimension 13 over Rationals>
> z:=LieCentre(g);
<two-sided ideal in <Lie algebra of dimension 13 over Rationals>, 
  (dimension 5)>
> Bz:=BasisVectors(Basis(z));
[ v.7+v.8+v.9+v.10+v.11+v.12+v.22, v.47+(-3)*v.48+(-1)*v.49+(-2)*v.50, 
  v.58+v.59, v.62, v.63 ]
> h*Bz[2];
(10)*v.47+(-30)*v.48+(-10)*v.49+(-20)*v.50
> h*Bz[3];
(14)*v.58+(14)*v.59
> h*Bz[4];
(16)*v.62
> h*Bz[5];
(18)*v.63
\end{verbatim}  
Les poids de {\tt z} sont 2,10,14,16,18; d'o{\`u} $m_{r}=18$. Il y a trois matrices {\`a}
{\'e}tudier.

\begin{enumerate}
\item $m_{i_2}=10$,  $m_{k(2)}=10$
\begin{verbatim}
> ((f*Bz[2])*Bz[2]);
(-150)*v.63
\end{verbatim}

\item $m_{i_3}=14$,  $m_{k(3)}=6$
\begin{verbatim}
> h*Bg[4];
(6)*v.30+(-12)*v.31+(6)*v.32+(6)*v.33+(-6)*v.35
> h*Bg[5];          
(6)*v.36+(-6)*v.37+(-6)*v.40

> ((f*Bg[4])*Bz[3]); 
(42)*v.63
\end{verbatim}

\item $m_{i_4}=16$,  $m_{k(4)}=4$
\begin{verbatim}
> h*Bg[3];
(4)*v.13+(4)*v.14+(4)*v.18+(-12)*v.26+(8)*v.28+(4)*v.29

> ((f*Bg[3])*Bz[4]);
(12)*v.63
\end{verbatim}
\end{enumerate}
{\bf Conclusion\,:} Cette orbite v{\'e}rifie la propri{\'e}t{\'e} $(P)$.\\\\

\item Caract{\'e}ristique\,:
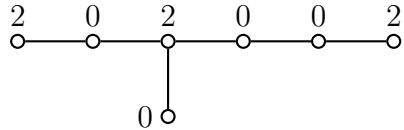
\begin{center}
\begin{pspicture}(-5,-1.2)(10,1)

\pscircle(-2,0){1mm}
\pscircle(-1,0){1mm}
\pscircle(0,0){1mm}
\pscircle(1,0){1mm}
\pscircle(2,0){1mm}
\pscircle(3,0){1mm}

\pscircle(0,-1){1mm}

\psline(-1.9,0)(-1.1,0)
\psline(-0.1,0)(-0.9,0)
\psline(0.1,0)(0.9,0)
\psline(1.1,0)(1.9,0)
\psline(2.1,0)(2.9,0)
\psline(0,-0.1)(0,-0.9)

\rput[b](-2,0.2){$2$}
\rput[b](-1,0.2){$0$}
\rput[b](0,0.2){$2$}
\rput[b](1,0.2){$0$}
\rput[b](2,0.2){$0$}
\rput[b](3,0.2){$2$}
\rput[r](-0.2,-1){$0$}
\end{pspicture}
\end{center} 
D{\'e}finition du $\mathfrak{sl}_{2}$-triplet\,:
\begin{verbatim}
> e:=x[8]+x[9]+x[13]+x[16]+x[17]+x[18]+x[29];                           
v.8+v.9+v.13+v.16+v.17+v.18+v.29
> f:=(14)*y[8]+(9)*y[9]+(-9)*y[10]+(11)*y[13]+(9)*y[16]
+(11)*y[17]+(8)*y[18]+(9)*y[19]+y[29];
(14)*v.71+(9)*v.72+(-9)*v.73+(11)*v.76+(9)*v.79
+(11)*v.80+(8)*v.81+(9)*v.82+v.92
> e*f;
(14)*v.127+(19)*v.128+(26)*v.129+(38)*v.130
+(29)*v.131+(20)*v.132+(11)*v.133
> h:=e*f;;
\end{verbatim}
Calcul de {\tt g}, {\tt Bg}, {\tt z} et {\tt Bz}\,:
\begin{verbatim}
> g:=LieCentralizer(L,Subspace(L,[e]));
<Lie algebra of dimension 17 over Rationals>
> z:=LieCentre(g);
<two-sided ideal in <Lie algebra of dimension 17 over Rationals>, 
  (dimension 3)>
> Bz:=BasisVectors(Basis(z)); 
[ v.8+v.9+v.13+v.16+v.17+v.18+v.29, v.56+(3)*v.57+(2)*v.59, v.63 ]
> h*Bz[2];
(10)*v.56+(30)*v.57+(20)*v.59
> h*Bz[3];
(14)*v.63
\end{verbatim}  
Les poids de {\tt z} sont 2,10,14; d'o{\`u} $m_{r}=14$. Il n'y a qu'une matrice {\`a} {\'e}tudier.

\begin{enumerate}
\item $m_{i_2}=10$,  $m_{k(2)}=6$
\begin{verbatim}
> h*Bg[7];
(6)*v.37+(-6)*v.38+(-6)*v.41
> h*Bg[8];
(6)*v.39+(6)*v.43+(-6)*v.45
> h*Bg[9];
(6)*v.42+(6)*v.46+(6)*v.49

> ((f*Bg[7])*Bz[2]);
(30)*v.63
\end{verbatim}
\end{enumerate}
{\bf Conclusion\,:} Cette orbite v{\'e}rifie la propri{\'e}t{\'e} $(P)$.\\

\item Caract{\'e}ristique\,:
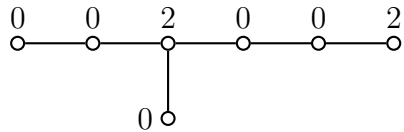
\begin{center}
\begin{pspicture}(-5,-1.2)(10,1)

\pscircle(-2,0){1mm}
\pscircle(-1,0){1mm}
\pscircle(0,0){1mm}
\pscircle(1,0){1mm}
\pscircle(2,0){1mm}
\pscircle(3,0){1mm}

\pscircle(0,-1){1mm}

\psline(-1.9,0)(-1.1,0)
\psline(-0.1,0)(-0.9,0)
\psline(0.1,0)(0.9,0)
\psline(1.1,0)(1.9,0)
\psline(2.1,0)(2.9,0)
\psline(0,-0.1)(0,-0.9)

\rput[b](-2,0.2){$0$}
\rput[b](-1,0.2){$0$}
\rput[b](0,0.2){$2$}
\rput[b](1,0.2){$0$}
\rput[b](2,0.2){$0$}
\rput[b](3,0.2){$2$}
\rput[r](-0.2,-1){$0$}
\end{pspicture}
\end{center} 
D{\'e}finition du $\mathfrak{sl}_{2}$-triplet\,:
\begin{verbatim}
> e:=x[13]+x[14]+x[15]+x[16]+x[17]+x[18]+x[33];                         
v.13+v.14+v.15+v.16+v.17+v.18+v.33
>f:=(9)*y[13]+(5)*y[14]+(2)*y[15]+(8)*y[16]+(8)*y[17]
+(2)*y[18]+(5)*y[33]; 
(9)*v.76+(5)*v.77+(2)*v.78+(8)*v.79+(8)*v.80+(2)*v.81+(5)*v.96
> e*f;
(10)*v.127+(15)*v.128+(20)*v.129+(30)*v.130
+(23)*v.131+(16)*v.132+(9)*v.133
> h:=e*f;;
\end{verbatim}
Calcul de {\tt g}, {\tt Bg}, {\tt z} et {\tt Bz}\,:
\begin{verbatim}
> g:=LieCentralizer(L,Subspace(L,[e]));
<Lie algebra of dimension 21 over Rationals>
> Bg:=BasisVectors(Basis(g));;
> z:=LieCentre(g);
<two-sided ideal in <Lie algebra of dimension 21 over Rationals>, 
  (dimension 4)>
> Bz:=BasisVectors(Basis(z));
[ v.13+v.14+v.15+v.16+v.17+v.18+v.33, v.61, v.62, v.63 ]
> h*Bz[2];
(10)*v.61
> h*Bz[3];
(10)*v.62
> h*Bz[4];
(10)*v.63
\end{verbatim}  
Les poids de {\tt z} sont 2,10,10,10; d'o{\`u} $m_{r}=10$. Il n'y a qu'une matrice {\`a}
{\'e}tudier.

\begin{enumerate}
\item $m_{i_2}=10$,  $m_{k(2)}=2$
\begin{verbatim}
> h*Bg[1];
(2)*v.15+(2)*v.18
> h*Bg[2]; 
(2)*v.4+(2/3)*v.19+(-4/3)*v.20+(2/3)*v.21+(2/3)*v.22+(-2/3)*v.23
> h*Bg[3];
(2)*v.19+(-4)*v.20+(2)*v.21+(2)*v.22+(4)*v.23+(-6)*v.24
> h*Bg[4];
(2)*v.7+(-6)*v.9+(4)*v.10+(2)*v.11+(2)*v.26+(-4)*v.27
> h*Bg[5];
(2)*v.9+(-2)*v.10+(2)*v.29
> h*Bg[6];
(2)*v.13+(2)*v.14+(2)*v.16+(2)*v.17+(2)*v.33

> ((f*Bg[1])*Bz[2]);((f*Bg[2])*Bz[2]);((f*Bg[3])*Bz[2]);
((f*Bg[4])*Bz[2]);((f*Bg[5])*Bz[2]);((f*Bg[6])*Bz[2]);
(-2)*v.61
(8/3)*v.62
(2)*v.62
(-4)*v.63
0*v.1
(-8)*v.61
>  ((f*Bg[1])*Bz[3]); ((f*Bg[2])*Bz[3]); ((f*Bg[3])*Bz[3]); 
((f*Bg[4])*Bz[3]);((f*Bg[5])*Bz[3]);((f*Bg[6])*Bz[3]);
(-2)*v.62
(4/3)*v.63
(4)*v.63
(2)*v.61
(2)*v.61
(-8)*v.62
> ((f*Bg[1])*Bz[4]); ((f*Bg[2])*Bz[4]); ((f*Bg[3])*Bz[4]); 
((f*Bg[4])*Bz[4]);((f*Bg[5])*Bz[4]);((f*Bg[6])*Bz[4]);
0*v.1
(-10/3)*v.61
(-10)*v.61
(10)*v.62
0*v.1
(-10)*v.63
\end{verbatim}
La matrice {\`a} {\'e}tudier est
$$\left[
\begin{array}{cccccc}
-2\alpha & -10/3 \gamma & -10 \gamma & 2 \beta & 2 \beta &-8 \alpha\\
-2\beta & -8/3 \alpha & 2\alpha & 10\gamma & 0 & -8\beta\\
0 & 4/3 \beta & 4\beta & -4\alpha & 0 & -10\gamma
\end{array}
\right]
 \cdot$$
Une {\'e}tude {\'e}l{\'e}mentaire montre que cette matrice est de rang 3 pour tout
 triplet $(\alpha,\beta,\gamma)$ non nul. 
\end{enumerate}
{\bf Conclusion\,:} Cette orbite v{\'e}rifie la propri{\'e}t{\'e} $(P)$.\\

\end{enumerate}
{\bf Conclusion\,:} Toutes les orbites nilpotentes distingu{\'e}es non
r{\'e}guli{\`e}res de $E_{7}$ v{\'e}rifient la propri{\'e}t{\'e} $(P)$.

\subsection{Calculs pour $E_{8}$}
\noindent D{\'e}finition de {\tt L}\,:
\begin{verbatim}
> L:=SimpleLieAlgebra("E",8,Rationals);
<Lie algebra of dimension 248 over Rationals>
> R:=RootSystem(L);
<root system of rank 7>
> P:=PositiveRoots(R);; 
> x:=PositiveRootVectors(R);
[ v.1, v.2, v.3, v.4, v.5, v.6, v.7, v.8, 
v.9, v.10, v.11, v.12, v.13, v.14, v.15, v.16, 
v.17, v.18, v.19, v.20, v.21, v.22, v.23, v.24, 
v.25, v.26, v.27, v.28, v.29, v.30, v.31, v.32, 
v.33, v.34, v.35, v.36, v.37, v.38, v.39, v.40, 
v.41, v.42, v.43, v.44, v.45, v.46, v.47, v.48, 
v.49, v.50, v.51, v.52, v.53, v.54, v.55, v.56, 
v.57, v.58, v.59, v.60, v.61, v.62, v.63, v.64, 
v.65, v.66, v.67, v.68, v.69, v.70, v.71, v.72, 
v.73, v.74, v.75, v.76, v.77, v.78, v.79, v.80, 
v.81, v.82, v.83, v.84, v.85, v.86, v.87, v.88, 
v.89, v.90, v.91, v.92, v.93, v.94, v.95, v.96, 
v.97, v.98, v.99, v.100, v.101, v.102, v.103, v.104, 
v.105, v.106, v.107, v.108, v.109, v.110, v.111, v.112, 
v.113, v.114, v.115, v.116, v.117, v.118, v.119, v.120 ]
> y:=NegativeRootVectors(R);
[ v.121, v.122, v.123, v.124, v.125, v.126, v.127, v.128, 
v.129, v.130, v.131, v.132, v.133, v.134, v.135, v.136, 
v.137, v.138, v.139, v.140, v.141, v.142, v.143, v.144, 
v.145, v.146, v.147, v.148, v.149, v.150, v.151, v.152, 
v.153, v.154, v.155, v.156, v.157, v.158, v.159, v.160, 
v.161, v.162, v.163, v.164, v.165, v.166, v.167, v.168, 
v.169, v.170, v.171, v.172, v.173, v.174, v.175, v.176, 
v.177, v.178, v.179, v.180, v.181, v.182, v.183, v.184, 
v.185, v.186, v.187, v.188, v.189, v.190, v.191, v.192, 
v.193, v.194, v.195, v.196, v.197, v.198, v.199, v.200, 
v.201, v.202, v.203, v.204, v.205, v.206, v.207, v.208, 
v.209, v.210, v.211, v.212, v.213, v.214, v.215, v.216, 
v.217, v.218, v.219, v.220, v.221, v.222, v.223, v.224, 
v.225, v.226, v.227, v.228, v.229, v.230, v.231, v.232, 
v.233, v.234, v.235, v.236, v.237, v.238, v.239, v.240 ]
> CanonicalGenerators(R)[3];
[ v.241, v.242, v.243, v.244, v.245, v.246, v.247, v.248 ]
\end{verbatim}
Dans $E_{8}$, il y a dix orbites nilpotentes distingu{\'e}es non
r{\'e}guli{\`e}res\,:\\

\begin{enumerate}
\item Caract{\'e}ristique\,:
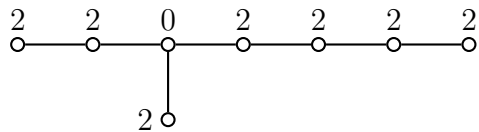
\begin{center}
\begin{pspicture}(-5,-1.2)(10,1)

\pscircle(-2,0){1mm}
\pscircle(-1,0){1mm}
\pscircle(0,0){1mm}
\pscircle(1,0){1mm}
\pscircle(2,0){1mm}
\pscircle(3,0){1mm}
\pscircle(4,0){1mm}

\pscircle(0,-1){1mm}

\psline(-1.9,0)(-1.1,0)
\psline(-0.1,0)(-0.9,0)
\psline(0.1,0)(0.9,0)
\psline(1.1,0)(1.9,0)
\psline(2.1,0)(2.9,0)
\psline(3.1,0)(3.9,0)
\psline(0,-0.1)(0,-0.9)

\rput[b](-2,0.2){$2$}
\rput[b](-1,0.2){$2$}
\rput[b](0,0.2){$0$}
\rput[b](1,0.2){$2$}
\rput[b](2,0.2){$2$}
\rput[b](3,0.2){$2$}
\rput[b](4,0.2){$2$}
\rput[r](-0.2,-1){$2$}
\end{pspicture}
\end{center}
D{\'e}finition du $\mathfrak{sl}_{2}$-triplet\,:
\begin{verbatim}
> e:=x[1]+x[2]+x[3]+x[6]+x[7]+x[8]+x[10]+x[12]; 
v.1+v.2+v.3+v.6+v.7+v.8+v.10+v.12
> f:=(72)*y[1]+(60)*y[2]+(142)*y[3]+(68)*y[5]+(132)*y[6]
+(90)*y[7]+(46)*y[8]+(38)*y[10]+(-38)*y[11]+(172)*y[12];                       
(72)*v.121+(68)*v.122+(142)*v.123+(68)*v.125+(132)*v.126
+(90)*v.127+(46)*v.128+(38)*v.130+(-38)*v.131+(172)*v.132
> e*f;
(72)*v.241+(106)*v.242+(142)*v.243+(210)*v.244
+(172)*v.245+(132)*v.246+(90)*v.247+(46)*v.248
> h:=e*f;;
\end{verbatim}
Calcul de {\tt g}, {\tt Bg}, {\tt z} et {\tt Bz}\,:
\begin{verbatim}
> g:=LieCentralizer(L,Subspace(L,[e]));
<Lie algebra of dimension 10 over Rationals>
> Bg:=BasisVectors(Basis(g));;
> z:=LieCentre(g);
<two-sided ideal in <Lie algebra of dimension 10 over Rationals>, 
  (dimension 7)>
> Bz:=BasisVectors(Basis(z));
[ v.1+v.2+v.3+v.6+v.7+v.8+v.10+v.12, v.54+(-1/2)*v.57+(-1/2)*v.58
+(-1/2)*v.59+(-1/2)*v.60+(-1/2)*v.61+(1/2)*v.62+(-1/2)*v.63, 
v.84+v.85+(-1)*v.86+(-1)*v.87+(2)*v.88, 
v.95+(-1)*v.96+v.97+(-1)*v.98+v.99,
v.109+(-1)*v.112, v.116, v.120 ]
> h*Bz[2];
(14)*v.54+(-7)*v.57+(-7)*v.58+(-7)*v.59+(-7)*v.60+(-7)*v.61
+(7)*v.62+(-7)*v.63
> h*Bz[3];
(22)*v.84+(22)*v.85+(-22)*v.86+(-22)*v.87+(44)*v.88
> h*Bz[4];
(26)*v.95+(-26)*v.96+(26)*v.97+(-26)*v.98+(26)*v.99
> h*Bz[5];
(34)*v.109+(-34)*v.112
> h*Bz[6];
(38)*v.116
> h*Bz[7];
(46)*v.120
\end{verbatim}  
Les poids de {\tt z} sont 2,14,22,26,34,38,46; d'o{\`u} $m_{r}=46$. Il y a cinq matrices {\`a}
{\'e}tudier.

\begin{enumerate}
\item $m_{i_2}=14$,  $m_{k(2)}=34$
\begin{verbatim}
> ((f*Bz[5])*Bz[2]);
(-119)*v.120
\end{verbatim}

\item $m_{i_3}=22$,  $m_{k(3)}=26$
\begin{verbatim}
> ((f*Bz[4])*Bz[3]);
(286)*v.120
\end{verbatim}

\item $m_{i_4}=26$,  $m_{k(4)}=22$
\begin{verbatim}
> ((f*Bz[3])*Bz[4]);
(286)*v.120
\end{verbatim}

\item $m_{i_5}=34$,  $m_{k(5)}=14$
\begin{verbatim}
> ((f*Bz[2])*Bz[5]);
(-119)*v.120
\end{verbatim}

\item $m_{i_6}=38$,  $m_{k(6)}=10$
\begin{verbatim}
> h*Bg[2];
(10)*v.38+(20)*v.39+(-10)*v.41+(-30)*v.42+(20)*v.43
+(40)*v.44+(40)*v.45+(10)*v.48+(20)*v.49

> ((f*Bg[2])*Bz[6]);               
(190)*v.120
\end{verbatim}
\end{enumerate}
{\bf Conclusion\,:} Cette orbite v{\'e}rifie la propri{\'e}t{\'e} $(P)$.\\

\item Caract{\'e}ristique\,:
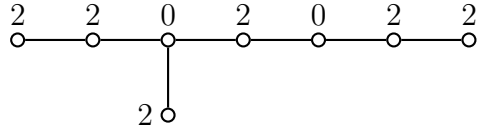
\begin{center}
\begin{pspicture}(-5,-1.2)(10,1)

\pscircle(-2,0){1mm}
\pscircle(-1,0){1mm}
\pscircle(0,0){1mm}
\pscircle(1,0){1mm}
\pscircle(2,0){1mm}
\pscircle(3,0){1mm}
\pscircle(4,0){1mm}

\pscircle(0,-1){1mm}

\psline(-1.9,0)(-1.1,0)
\psline(-0.1,0)(-0.9,0)
\psline(0.1,0)(0.9,0)
\psline(1.1,0)(1.9,0)
\psline(2.1,0)(2.9,0)
\psline(3.1,0)(3.9,0)
\psline(0,-0.1)(0,-0.9)

\rput[b](-2,0.2){$2$}
\rput[b](-1,0.2){$2$}
\rput[b](0,0.2){$0$}
\rput[b](1,0.2){$2$}
\rput[b](2,0.2){$0$}
\rput[b](3,0.2){$2$}
\rput[b](4,0.2){$2$}
\rput[r](-0.2,-1){$2$}
\end{pspicture}
\end{center}
D{\'e}finition du $\mathfrak{sl}_{2}$-triplet\,:
\begin{verbatim}
> e:=x[1]+x[2]+x[3]+x[5]+x[7]+x[8]+x[10]+x[20];                         
v.1+v.2+v.3+v.5+v.7+v.8+v.10+v.20
> f:=(61)*y[1]+(22)*y[2]+(118)*y[3]+(34)*y[5]+(74)*y[7]+(38)*y[8]
+(66)*y[10]+(-66)*y[11]+(22)*y[13]+(22)*y[14]+(108)*y[20];            
(60)*v.121+(22)*v.122+(118)*v.123+(34)*v.125+(74)*v.127+
(38)*v.128+(66)*v.130+(-66)*v.131+(22)*v.133+(22)*v.134+(108)*v.140
> e*f;
(60)*v.241+(88)*v.242+(118)*v.243+(174)*v.244+(142)*v.245+(108)*v.246
+(74)*v.247+(38)*v.248
> h:=e*f;;
\end{verbatim}
Calcul de {\tt g}, {\tt Bg}, {\tt z} et {\tt Bz}\,:
\begin{verbatim}
> g:=LieCentralizer(L,Subspace(L,[e]));
<Lie algebra of dimension 12 over Rationals>
> Bg:=BasisVectors(Basis(g));;
> z:=LieCentre(g);
<two-sided ideal in <Lie algebra of dimension 12 over Rationals>, 
  (dimension 6)>
> Bz:=BasisVectors(Basis(z));
[ v.1+v.2+v.3+v.5+v.7+v.8+v.10+v.20, 
v.64+v.65+(2)*v.67+v.69+v.71+v.73+v.74+v.76, 
v.97+(-1)*v.98+v.99+v.100, 
v.104+v.107+(-1)*v.108+(-1)*v.110, v.117, v.120 ]
> h*Bz[2];
(14)*v.64+(14)*v.65+(28)*v.67+(14)*v.69+(14)*v.71+
(14)*v.73+(14)*v.74+(14)*v.76
> h*Bz[3];
(22)*v.97+(-22)*v.98+(22)*v.99+(22)*v.100 
> h*Bz[4];
(26)*v.104+(26)*v.107+(-26)*v.108+(-26)*v.110
> h*Bz[5];
(34)*v.117
> h*Bz[6];
(38)*v.120
\end{verbatim}  
Les poids de {\tt z} sont 2,14,22,26,34,38; d'o{\`u} $m_{r}=38$. Il y a quatre matrices {\`a}
{\'e}tudier.

\begin{enumerate}
\item $m_{i_2}=14$,  $m_{k(2)}=26$
\begin{verbatim}
> ((f*Bz[4])*Bz[2]);
(-182)*v.120
\end{verbatim}

\item $m_{i_3}=22$,  $m_{k(3)}=18$
\begin{verbatim}
> h*Bg[6];
(18)*v.81+(36)*v.85+(-18)*v.86+(-18)*v.87+(18)*v.88

> ((f*Bg[6])*Bz[3]);                                                    
(-198)*v.120
\end{verbatim}

\item $m_{i_4}=26$,  $m_{k(4)}=14$
\begin{verbatim}
> ((f*Bz[2])*Bz[4]);
(-182)*v.120
\end{verbatim}

\item $m_{i_5}=34$,  $m_{k(5)}=6$
\begin{verbatim}
> h*Bg[2];
(6)*v.23+(9)*v.24+(-3)*v.25+(9)*v.29+(-3)*v.31
+(-9)*v.32+(3)*v.33+(6)*v.34
+(-3)*v.35+(-3)*v.36+(-3)*v.40

> ((f*Bg[2])*Bz[5]);
(51)*v.120
\end{verbatim}
\end{enumerate}
{\bf Conclusion\,:} Cette orbite v{\'e}rifie la propri{\'e}t{\'e} $(P)$.\\

\item Caract{\'e}ristique\,:
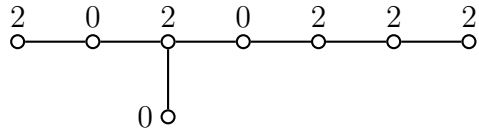
\begin{center}
\begin{pspicture}(-5,-1.2)(10,1)

\pscircle(-2,0){1mm}
\pscircle(-1,0){1mm}
\pscircle(0,0){1mm}
\pscircle(1,0){1mm}
\pscircle(2,0){1mm}
\pscircle(3,0){1mm}
\pscircle(4,0){1mm}

\pscircle(0,-1){1mm}

\psline(-1.9,0)(-1.1,0)
\psline(-0.1,0)(-0.9,0)
\psline(0.1,0)(0.9,0)
\psline(1.1,0)(1.9,0)
\psline(2.1,0)(2.9,0)
\psline(3.1,0)(3.9,0)
\psline(0,-0.1)(0,-0.9)

\rput[b](-2,0.2){$2$}
\rput[b](-1,0.2){$0$}
\rput[b](0,0.2){$2$}
\rput[b](1,0.2){$0$}
\rput[b](2,0.2){$2$}
\rput[b](3,0.2){$2$}
\rput[b](4,0.2){$2$}
\rput[r](-0.2,-1){$0$}
\end{pspicture}
\end{center}
D{\'e}finition du $\mathfrak{sl}_{2}$-triplet\,:
\begin{verbatim}
> e:=x[7]+x[8]+x[9]+x[10]+x[11]+x[12]+x[13]+x[25];                      
v.7+v.8+v.9+v.10+v.11+v.12+v.13+v.25
> f:=(66)*y[7]+(34)*y[8]+(52)*y[9]+(75)*y[10]+(49)*y[11]+(27)*y[12]
+(96)*y[13]+y[25];                                                  
(66)*v.127+(34)*v.128+(52)*v.129+(75)*v.130+(49)*v.131+(27)*v.132
+(96)*v.133+v.145
> e*f;
(52)*v.241+(76)*v.242+(102)*v.243+(152)*v.244
+(124)*v.245+(96)*v.246+(66)*v.247+(34)*v.248
> h:=e*f;;
\end{verbatim}
Calcul de {\tt g}, {\tt Bg}, {\tt z} et {\tt Bz}\,:
\begin{verbatim}
> g:=LieCentralizer(L,Subspace(L,[e]));
<Lie algebra of dimension 14 over Rationals>
> Bg:=BasisVectors(Basis(g));;
> z:=LieCentre(g);
<two-sided ideal in <Lie algebra of dimension 14 over Rationals>, 
  (dimension 6)>
> Bz:=BasisVectors(Basis(z));
[ v.7+v.8+v.9+v.10+v.11+v.12+v.13+v.25,
v.74+(-2)*v.77+(-1)*v.78+(-1)*v.80
+(-1)*v.82, v.104+(-1)*v.105+(-1)*v.106, v.113+v.114, v.117, v.120 ]
> h*Bz[2];
(14)*v.74+(-28)*v.77+(-14)*v.78+(-14)*v.80+(-14)*v.82
> h*Bz[3];
(22)*v.104+(-22)*v.105+(-22)*v.106
> h*Bz[4];
(26)*v.113+(26)*v.114
> h*Bz[5];
(28)*v.117
> h*Bz[6];
(34)*v.120
\end{verbatim}  
Les poids de {\tt z} sont 2,14,22,26,28,34; d'o{\`u} $m_{r}=34$. Il y a quatre matrices {\`a}
{\'e}tudier.

\begin{enumerate}
\item $m_{i_2}=14$,  $m_{k(2)}=22$
\begin{verbatim}
> ((f*Bz[3])*Bz[2]);
(-154)*v.120
\end{verbatim}

\item $m_{i_3}=22$,  $m_{k(3)}=14$
\begin{verbatim}
> ((f*Bz[2])*Bz[3]);
(-154)*v.120
\end{verbatim}

\item $m_{i_4}=26$,  $m_{k(4)}=10$
\begin{verbatim}
> h*Bg[4];
(10)*v.54+(-5)*v.56+(-5)*v.58+(-10)*v.59+(5)*v.61+(-15)*v.63
> h*Bg[5];
(10)*v.62+(-10)*v.64+(-10)*v.69

> ((f*Bg[4])*Bz[4]);
(-65)*v.120
\end{verbatim}

\item $m_{i_5}=28$,  $m_{k(5)}=8$
\begin{verbatim}
> h*Bg[3];
(8)*v.36+(8)*v.39+(8)*v.45+(-24)*v.50+(16)*v.55+(-8)*v.57

> ((f*Bg[3])*Bz[5]);
(24)*v.120
\end{verbatim}
\end{enumerate}
{\bf Conclusion\,:} Cette orbite v{\'e}rifie la propri{\'e}t{\'e} $(P)$.\\

\item Caract{\'e}ristique\,:
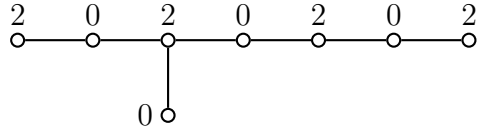
\begin{center}
\begin{pspicture}(-5,-1.2)(10,1)

\pscircle(-2,0){1mm}
\pscircle(-1,0){1mm}
\pscircle(0,0){1mm}
\pscircle(1,0){1mm}
\pscircle(2,0){1mm}
\pscircle(3,0){1mm}
\pscircle(4,0){1mm}

\pscircle(0,-1){1mm}

\psline(-1.9,0)(-1.1,0)
\psline(-0.1,0)(-0.9,0)
\psline(0.1,0)(0.9,0)
\psline(1.1,0)(1.9,0)
\psline(2.1,0)(2.9,0)
\psline(3.1,0)(3.9,0)
\psline(0,-0.1)(0,-0.9)

\rput[b](-2,0.2){$2$}
\rput[b](-1,0.2){$0$}
\rput[b](0,0.2){$2$}
\rput[b](1,0.2){$0$}
\rput[b](2,0.2){$2$}
\rput[b](3,0.2){$0$}
\rput[b](4,0.2){$2$}
\rput[r](-0.2,-1){$0$}
\end{pspicture}
\end{center}
D{\'e}finition du $\mathfrak{sl}_{2}$-triplet\,:
\begin{verbatim}
> e:=x[9]+x[10]+x[11]+x[12]+x[13]+x[14]+x[15]+x[25];                    
v.9+v.10+v.11+v.12+v.13+v.14+v.15+v.25
> f:=(44)*y[9]+(50)*y[10]+(28)*y[11]+(36)*y[12]+(54)*y[13]+(26)*y[14]
+(28)*y[15]+(14)*y[25];                                           
(44)*v.129+(50)*v.130+(28)*v.131+(36)*v.132+(54)*v.133+(26)*v.134
+(28)*v.135+(14)*v.145
> e*f;
(44)*v.241+(64)*v.242+(86)*v.243+(128)*v.244+(104)*v.245+(80)*v.246
+(54)*v.247+(28)*v.248
> h:=e*f;;
\end{verbatim}
Calcul de {\tt g}, {\tt Bg}, {\tt z} et {\tt Bz}\,:
\begin{verbatim}
> g:=LieCentralizer(L,Subspace(L,[e]));
<Lie algebra of dimension 16 over Rationals>
> Bg:=BasisVectors(Basis(g));;
> z:=LieCentre(g);
<two-sided ideal in <Lie algebra of dimension 16 over Rationals>, 
  (dimension 5)>
> Bz:=BasisVectors(Basis(z));
[ v.9+v.10+v.11+v.12+v.13+v.14+v.15+v.25, 
v.86+v.87+v.88+(4)*v.93+(-3)*v.95, 
v.111+v.115, v.119, v.120 ]
> h*Bz[2];
(14)*v.86+(14)*v.87+(14)*v.88+(56)*v.93+(-42)*v.95
> h*Bz[3];
(22)*v.111+(22)*v.115
> h*Bz[4];
(26)*v.119
> h*Bz[5];
(28)*v.120
\end{verbatim}  
Les poids de {\tt z} sont 2,14,22,26,28; d'o{\`u} $m_{r}=28$. Il y a trois matrices {\`a}
{\'e}tudier.

\begin{enumerate}
\item $m_{i_2}=14$,  $m_{k(2)}=16$
\begin{verbatim}
> h*Bg[10];
(16)*v.97+(16)*v.98+(-16)*v.99+(-16)*v.100

> ((f*Bg[10])*Bz[2]);
(-112)*v.120
\end{verbatim}

\item $m_{i_3}=22$,  $m_{k(3)}=8$
\begin{verbatim}
> h*Bg[4];
(8)*v.45+(8)*v.47+(-8)*v.57+(-24)*v.58+(24)*v.59+(16)*v.61+(16)*v.62
> h*Bg[5];

> ((f*Bg[4])*Bz[3]); 
(88)*v.120
\end{verbatim}

\item $m_{i_4}=26$,  $m_{k(4)}=4$
\begin{verbatim}
> h*Bg[2];
(4)*v.16+(4)*v.20+(4)*v.22+(-12)*v.30
+(8)*v.32+(4)*v.33+(-12)*v.34+(20)*v.35

> ((f*Bg[2])*Bz[4]);
(26)*v.120
\end{verbatim}
\end{enumerate}
{\bf Conclusion\,:} Cette orbite v{\'e}rifie la propri{\'e}t{\'e} $(P)$.\\

\item Caract{\'e}ristique\,:
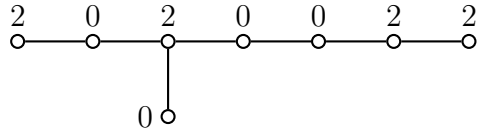
\begin{center}
\begin{pspicture}(-5,-1.2)(10,1)

\pscircle(-2,0){1mm}
\pscircle(-1,0){1mm}
\pscircle(0,0){1mm}
\pscircle(1,0){1mm}
\pscircle(2,0){1mm}
\pscircle(3,0){1mm}
\pscircle(4,0){1mm}

\pscircle(0,-1){1mm}

\psline(-1.9,0)(-1.1,0)
\psline(-0.1,0)(-0.9,0)
\psline(0.1,0)(0.9,0)
\psline(1.1,0)(1.9,0)
\psline(2.1,0)(2.9,0)
\psline(3.1,0)(3.9,0)
\psline(0,-0.1)(0,-0.9)

\rput[b](-2,0.2){$2$}
\rput[b](-1,0.2){$0$}
\rput[b](0,0.2){$2$}
\rput[b](1,0.2){$0$}
\rput[b](2,0.2){$0$}
\rput[b](3,0.2){$2$}
\rput[b](4,0.2){$2$}
\rput[r](-0.2,-1){$0$}
\end{pspicture}
\end{center}
D{\'e}finition du $\mathfrak{sl}_{2}$-triplet\,:
\begin{verbatim}
> e:=x[8]+x[9]+x[10]+x[14]+x[18]+x[19]+x[20]+x[33];                     
v.8+v.9+v.10+v.14+v.18+v.19+v.20+v.33
> f:=(8)*y[8]+(40)*y[9]+(22)*y[10]+(-22)*y[11]+(50)*y[14]+(35)*y[18]
+(37)*y[19]+(21)*y[20]+(35)*y[21]+y[33];                          
(26)*v.128+(40)*v.129+(22)*v.130+(-22)*v.131+(50)*v.134+(37)*v.138
+(37)*v.139+(21)*v.140+(35)*v.141+v.153
> e*f;                                                                  
(40)*v.241+(58)*v.242+(78)*v.243+(116)*v.244+(94)*v.245+(72)*v.246
+(50)*v.247+(26)*v.248
> h:=e*f;;
\end{verbatim}
Calcul de {\tt g}, {\tt Bg}, {\tt z} et {\tt Bz}\,:
\begin{verbatim}
> g:=LieCentralizer(L,Subspace(L,[e]));
<Lie algebra of dimension 18 over Rationals>
> Bg:=BasisVectors(Basis(g));;
> z:=LieCentre(g);
<two-sided ideal in <Lie algebra of dimension 18 over Rationals>, 
  (dimension 4)>
> Bz:=BasisVectors(Basis(z));
[ v.8+v.9+v.10+v.14+v.18+v.19+v.20+v.33, 
v.94+v.95+v.96+v.97, v.116, v.120 ]
> h*Bz[2];
(14)*v.94+(14)*v.95+(14)*v.96+(14)*v.97
> h*Bz[3];
(22)*v.116
> h*Bz[4];
(26)*v.120
\end{verbatim}  
Les poids de {\tt z} sont 2,14,22,26; d'o{\`u} $m_{r}=26$. Il y a deux matrices {\`a} {\'e}tudier.

\begin{enumerate}
\item $m_{i_2}=14$,  $m_{k(2)}=14$
\begin{verbatim}
> ((f*Bz[2])*Bz[2]);
(-98)*v.120
\end{verbatim}

\item $m_{i_3}=22$,  $m_{k(3)}=6$
\begin{verbatim}
> h*Bg[4];
(6)*v.42+(-12)*v.43+(12)*v.44+(-12)*v.45
+(6)*v.46+(-12)*v.49+(6)*v.52
+(-6)*v.55
> h*Bg[5];
(6)*v.50+(-6)*v.51+(-6)*v.57+(-6)*v.61

> ((f*Bg[4])*Bz[3]);                
(-66)*v.120
\end{verbatim}
\end{enumerate}
{\bf Conclusion\,:} Cette orbite v{\'e}rifie la propri{\'e}t{\'e} $(P)$.\\

\item Caract{\'e}ristique\,:
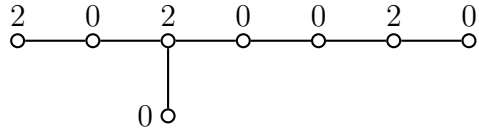
\begin{center}
\begin{pspicture}(-5,-1.2)(10,1)

\pscircle(-2,0){1mm}
\pscircle(-1,0){1mm}
\pscircle(0,0){1mm}
\pscircle(1,0){1mm}
\pscircle(2,0){1mm}
\pscircle(3,0){1mm}
\pscircle(4,0){1mm}

\pscircle(0,-1){1mm}

\psline(-1.9,0)(-1.1,0)
\psline(-0.1,0)(-0.9,0)
\psline(0.1,0)(0.9,0)
\psline(1.1,0)(1.9,0)
\psline(2.1,0)(2.9,0)
\psline(3.1,0)(3.9,0)
\psline(0,-0.1)(0,-0.9)

\rput[b](-2,0.2){$2$}
\rput[b](-1,0.2){$0$}
\rput[b](0,0.2){$2$}
\rput[b](1,0.2){$0$}
\rput[b](2,0.2){$0$}
\rput[b](3,0.2){$2$}
\rput[b](4,0.2){$0$}
\rput[r](-0.2,-1){$0$}
\end{pspicture}
\end{center}
D{\'e}finition du $\mathfrak{sl}_{2}$-triplet\,:
\begin{verbatim}
> e:=x[9]+x[10]+x[14]+x[15]+x[18]+x[19]+x[20]+x[33];                    
v.9+v.10+v.14+v.15+v.18+v.19+v.20+v.33
> f:=(36)*y[9]+(20)*y[10]+(-20)*y[11]+(22)*y[14]+(22)*y[15]+(20)*y[18]
+(22)*y[19]+(30)*y[20]+(20)*y[21]+(12)*y[33];                    
(36)*v.129+(20)*v.130+(-20)*v.131+(22)*v.134+(22)*v.135+(20)*v.138
+(22)*v.139+(30)*v.140+(20)*v.141+(12)*v.153
> e*f;
(36)*v.241+(52)*v.242+(70)*v.243+(104)*v.244+(84)*v.245+(64)*v.246
+(44)*v.247+(22)*v.248
> h:=e*f;;
\end{verbatim}
Calcul de {\tt g}, {\tt Bg}, {\tt z} et {\tt Bz}\,:
\begin{verbatim}
> g:=LieCentralizer(L,Subspace(L,[e]));
<Lie algebra of dimension 20 over Rationals>
> Bg:=BasisVectors(Basis(g));;
> z:=LieCentre(g);
<two-sided ideal in <Lie algebra of dimension 20 over Rationals>, 
  (dimension 4)>
> Bz:=BasisVectors(Basis(z));
[ v.9+v.10+v.14+v.15+v.18+v.19+v.20+v.33, 
v.97+(-1)*v.99+(3)*v.105+(4)*v.108,
v.119, v.120 ]
> h*Bz[2];                   
(14)*v.97+(-14)*v.99+(42)*v.105+(56)*v.108
> h*Bz[3];
(22)*v.119
> h*Bz[4];
(22)*v.120
\end{verbatim}  
Les poids de {\tt z} sont 2,14,22,22; d'o{\`u} $m_{r}=22$. Il y a deux matrices {\`a} {\'e}tudier.

\begin{enumerate}
\item $m_{i_2}=14$,  $m_{k(2)}=10$
\begin{verbatim}
> h*Bg[7];
(10)*v.70+(10)*v.80+(10)*v.81+(10)*v.83
> h*Bg[8];
(10)*v.75+(10)*v.76+(10)*v.77+(10)*v.86
> h*Bg[9];
(10)*v.81+(20)*v.83+(10)*v.84+(-10)*v.85+(10)*v.87
> h*Bg[10];
(10)*v.76+(10)*v.77+(-10)*v.79+(10)*v.82+(10)*v.86+(-10)*v.90

> ((f*Bg[7])*Bz[2]);
0*v.1
> ((f*Bg[8])*Bz[2]);
0*v.1
> ((f*Bg[9])*Bz[2]);
(-70)*v.119
> ((f*Bg[10])*Bz[2]); 
(70)*v.120
\end{verbatim}
La matrice correspondante est 
$$\left[
\begin{array}{cccc}
0 & 0 & -70 & 0 \\
0 & 0 & 0 & 70
\end{array}
\right] \cdot$$
C'est clairement une matrice de rang 2.

\item $m_{i_3}=22$,  $m_{k(3)}=2$
\begin{verbatim}
> h*Bg[3];
(2)*v.9+(2)*v.10+(2)*v.14+(2)*v.15
+(2)*v.18+(2)*v.19+(2)*v.20+(2)*v.33

> ((f*Bg[1])*Bz[3]); ((f*Bg[2])*Bz[3]);((f*Bg[3])*Bz[3]);
(-2)*v.119
(14)*v.120
(-22)*v.119
> ((f*Bg[1])*Bz[4]); ((f*Bg[2])*Bz[4]); ((f*Bg[3])*Bz[4]);
0*v.1
(22)*v.119
(-22)*v.120
\end{verbatim}
La matrice {\`a} {\'e}tudier est
$$\left[
\begin{array}{cccc}
-2\alpha & 22 \beta &-22 \alpha\\
0 & 14\alpha & -22 \beta
\end{array}
\right]
 \cdot$$
Cette matrice est de rang 2 pour tout couple
 $(\alpha,\beta)$ non nul. 
\end{enumerate}
{\bf Conclusion\,:} Cette orbite v{\'e}rifie la propri{\'e}t{\'e} $(P)$.\\

\item Caract{\'e}ristique\,:
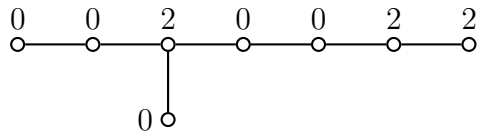
\begin{center}
\begin{pspicture}(-5,-1.2)(10,1)

\pscircle(-2,0){1mm}
\pscircle(-1,0){1mm}
\pscircle(0,0){1mm}
\pscircle(1,0){1mm}
\pscircle(2,0){1mm}
\pscircle(3,0){1mm}
\pscircle(4,0){1mm}

\pscircle(0,-1){1mm}

\psline(-1.9,0)(-1.1,0)
\psline(-0.1,0)(-0.9,0)
\psline(0.1,0)(0.9,0)
\psline(1.1,0)(1.9,0)
\psline(2.1,0)(2.9,0)
\psline(3.1,0)(3.9,0)
\psline(0,-0.1)(0,-0.9)

\rput[b](-2,0.2){$0$}
\rput[b](-1,0.2){$0$}
\rput[b](0,0.2){$2$}
\rput[b](1,0.2){$0$}
\rput[b](2,0.2){$0$}
\rput[b](3,0.2){$2$}
\rput[b](4,0.2){$2$}
\rput[r](-0.2,-1){$0$}
\end{pspicture}
\end{center}
D{\'e}finition du $\mathfrak{sl}_{2}$-triplet\,:
\begin{verbatim}
> e:=x[8]+x[14]+x[16]+x[17]+x[18]+x[19]+x[20]+x[38];
v.8+v.14+v.16+v.17+v.18+v.19+v.20+v.38
> f:=(22)*y[8]+(42)*y[14]+(16)*y[16]+(2)*y[17]+(30)*y[18]+(30)*y[19]
+(2)*y[20]+(16)*y[38];                                                   
(22)*v.128+(42)*v.134+(16)*v.136+(2)*v.137+(30)*v.138+(30)*v.139
+(2)*v.140+(16)*v.158
> e*f;
(32)*v.241+(48)*v.242+(64)*v.243+(96)*v.244+(78)*v.245+(60)*v.246
+(42)*v.247+(22)*v.248
> h:=e*f;;
\end{verbatim}
Calcul de {\tt g}, {\tt Bg}, {\tt z} et {\tt Bz}\,:
\begin{verbatim}
> g:=LieCentralizer(L,Subspace(L,[e]));
<Lie algebra of dimension 22 over Rationals>
> Bg:=BasisVectors(Basis(g));;
> z:=LieCentre(g);
<two-sided ideal in <Lie algebra of dimension 22 over Rationals>, 
  (dimension 5)>
> Bz:=BasisVectors(Basis(z));
[ v.8+v.14+v.16+v.17+v.18+v.19+v.20+v.38, 
v.99+(-1)*v.107+(-1)*v.108, v.116, 
v.118, v.120 ]
> h*Bz[2];
(14)*v.99+(-14)*v.107+(-14)*v.108
> h*Bz[3];
(18)*v.116
> h*Bz[4];
(18)*v.118
> h*Bz[5];
(22)*v.120
\end{verbatim}  
Les poids de {\tt z} sont 2,14,18,18,22; d'o{\`u} $m_{r}=22$. Il y a deux matrices {\`a}
{\'e}tudier.

\begin{enumerate}
\item $m_{i_2}=14$,  $m_{k(2)}=10$
\begin{verbatim}
> h*Bg[11];
(10)*v.83+(10)*v.89
> h*Bg[12];
(10)*v.86+(-10)*v.93
> h*Bg[13];
(10)*v.74+(-10)*v.77+(-10)*v.90+(20)*v.97

> ((f*Bg[13])*Bz[2]);
(-70)*v.120
\end{verbatim}

\item $m_{i_3}=18$,  $m_{k(3)}=6$
\begin{verbatim}
> h*Bg[6];
(6)*v.47+(-12)*v.50+(-6)*v.58+(6)*v.59+(12)*v.61+(18)*v.63
> h*Bg[7];
(6)*v.36+(-6)*v.49+(-3)*v.54+(3)*v.64+(-3)*v.66+(-9)*v.69

> ((f*Bg[6])*Bz[3]); ((f*Bg[7])*Bz[3]);
(-12)*v.120
0*v.1
> ((f*Bg[6])*Bz[4]); ((f*Bg[7])*Bz[4]);
0*v.1
(6)*v.120
\end{verbatim}
La matrice {\`a} {\'e}tudier est 
$\left[
\begin{array}{cc}
-12 \alpha & 6\beta 
\end{array}
\right]$; elle est de rang 1 si le couple $(\alpha,\beta)$ est non
nul. 
\end{enumerate}
{\bf Conclusion\,:} Cette orbite v{\'e}rifie la propri{\'e}t{\'e} $(P)$.\\

\item Caract{\'e}ristique\,:
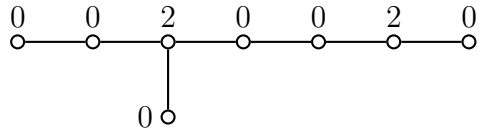
\begin{center}
\begin{pspicture}(-5,-1.2)(10,1)

\pscircle(-2,0){1mm}
\pscircle(-1,0){1mm}
\pscircle(0,0){1mm}
\pscircle(1,0){1mm}
\pscircle(2,0){1mm}
\pscircle(3,0){1mm}
\pscircle(4,0){1mm}

\pscircle(0,-1){1mm}

\psline(-1.9,0)(-1.1,0)
\psline(-0.1,0)(-0.9,0)
\psline(0.1,0)(0.9,0)
\psline(1.1,0)(1.9,0)
\psline(2.1,0)(2.9,0)
\psline(3.1,0)(3.9,0)
\psline(0,-0.1)(0,-0.9)

\rput[b](-2,0.2){$0$}
\rput[b](-1,0.2){$0$}
\rput[b](0,0.2){$2$}
\rput[b](1,0.2){$0$}
\rput[b](2,0.2){$0$}
\rput[b](3,0.2){$2$}
\rput[b](4,0.2){$0$}
\rput[r](-0.2,-1){$0$}
\end{pspicture}
\end{center}
D{\'e}finition du $\mathfrak{sl}_{2}$-triplet\,:
\begin{verbatim}
> e:=x[14]+x[15]+x[16]+x[17]+x[18]+x[19]+x[20]+x[38];                      
v.14+v.15+v.16+v.17+v.18+v.19+v.20+v.38
> f:=(18)*y[14]+(18)*y[15]+(8)*y[16]+(8)*y[17]+(14)*y[18]+(20)*y[19]
+(14)*y[20]+(20)*y[38];                                                  
(18)*v.134+(18)*v.135+(8)*v.136+(8)*v.137
+(14)*v.138+(20)*v.139+(14)*v.140+(20)*v.158
> e*f;
(28)*v.241+(42)*v.242+(56)*v.243+(84)*v.244
+(68)*v.245+(52)*v.246+(36)*v.247+(18)*v.248
> h:=e*f;;
\end{verbatim}
Calcul de {\tt g}, {\tt Bg}, {\tt z} et {\tt Bz}\,:
\begin{verbatim}
> g:=LieCentralizer(L,Subspace(L,[e]));
<Lie algebra of dimension 24 over Rationals>
> Bg:=BasisVectors(Basis(g));;
> z:=LieCentre(g);
<two-sided ideal in <Lie algebra of dimension 24 over Rationals>, 
  (dimension 4)>
> Bz:=BasisVectors(Basis(z));
[ v.14+v.15+v.16+v.17+v.18+v.19+v.20+v.38, v.111+v.112, v.119, v.120 ]
> h*Bz[2];                   
(14)*v.111+(14)*v.112
> h*Bz[3];
(18)*v.119
> h*Bz[4];
(18)*v.120
\end{verbatim}  
Les poids de {\tt z} sont 2,14,18,18; d'o{\`u} $m_{r}=18$. Il y a deux matrices {\`a} {\'e}tudier.

\begin{enumerate}
\item $m_{i_2}=14$,  $m_{k(2)}=6$
\begin{verbatim}
> h*Bg[5];
(6)*v.58+(-6)*v.59+(6)*v.60+(-6)*v.63
> h*Bg[6];
(6)*v.49+(6)*v.64+(6)*v.65+(-6)*v.66+(-6)*v.67+(6)*v.68
> h*Bg[7];
(6)*v.64+(12)*v.65+(-6)*v.66+(-6)*v.67+(12)*v.68+(-6)*v.69
> h*Bg[8];
(6)*v.53+(-6)*v.55+(-6)*v.56+(-6)*v.71+(-6)*v.72+(6)*v.73
> h*Bg[9];
(6)*v.60+(-3)*v.61+(-3)*v.62+(-3)*v.63+(-3)*v.78

> ((f*Bg[5])*Bz[2]);
(14)*v.119
> ((f*Bg[6])*Bz[2]);
0*v.1
> ((f*Bg[7])*Bz[2]);
(14)*v.120
> ((f*Bg[8])*Bz[2]);
0*v.1
> ((f*Bg[9])*Bz[2]);
(7)*v.119
\end{verbatim}
 matrice {\`a} {\'e}tudier est 
$$\left[
\begin{array}{cccccc}
14 & 0 & 0 & 0 & 7 \\
0 & 0 & 14  & 0
\end{array}
\right] \cdot $$ 
C'est clairement une matrice de rang 2.

\item $m_{i_3}=18$,  $m_{k(3)}=2$
\begin{verbatim}
> h*Bg[2];
(2)*v.7+(-10)*v.10+(8)*v.11+(2)*v.12
+(-4)*v.29+(2)*v.30+(-4)*v.31+(-4)*v.33
> h*Bg[3];
(2)*v.14+(2)*v.15+(2)*v.16+(2)*v.17
+(2)*v.18+(2)*v.19+(2)*v.20+(2)*v.38

> ((f*Bg[1])*Bz[3]); ((f*Bg[2])*Bz[3]);((f*Bg[3])*bz[3]);
(-18/5)*v.120
0*v.1
(-18)*v.119
> ((f*Bg[1])*Bz[4]) ;((f*Bg[2])*Bz[4]); ((f*Bg[3])*Bz[4]);
0*v.1
(18)*v.119
(-18)*v.120
\end{verbatim}
La matrice {\`a} {\'e}tudier est 
$$\left[
\begin{array}{ccc}
0 & 18\beta & -18\alpha \\
-18/5 \alpha & 0 & -18\beta
\end{array}
\right] \cdot $$ 
C'est une matrice de rang 2 pour tout couple $(\alpha,\beta)$ non
nul. 
\end{enumerate}
{\bf Conclusion\,:} Cette orbite v{\'e}rifie la propri{\'e}t{\'e} $(P)$.\\

\item Caract{\'e}ristique\,:
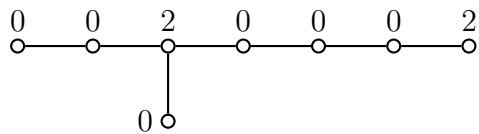
\begin{center}
\begin{pspicture}(-5,-1.2)(10,1)

\pscircle(-2,0){1mm}
\pscircle(-1,0){1mm}
\pscircle(0,0){1mm}
\pscircle(1,0){1mm}
\pscircle(2,0){1mm}
\pscircle(3,0){1mm}
\pscircle(4,0){1mm}

\pscircle(0,-1){1mm}

\psline(-1.9,0)(-1.1,0)
\psline(-0.1,0)(-0.9,0)
\psline(0.1,0)(0.9,0)
\psline(1.1,0)(1.9,0)
\psline(2.1,0)(2.9,0)
\psline(3.1,0)(3.9,0)
\psline(0,-0.1)(0,-0.9)

\rput[b](-2,0.2){$0$}
\rput[b](-1,0.2){$0$}
\rput[b](0,0.2){$2$}
\rput[b](1,0.2){$0$}
\rput[b](2,0.2){$0$}
\rput[b](3,0.2){$0$}
\rput[b](4,0.2){$2$}
\rput[r](-0.2,-1){$0$}
\end{pspicture}
\end{center}
D{\'e}finition du $\mathfrak{sl}_{2}$-triplet\,:
\begin{verbatim}
> e:=x[15]+x[16]+x[17]+x[18]+x[19]+x[20]+x[38]+x[46];                      
v.15+v.16+v.17+v.18+v.19+v.20+v.38+v.46
> f:=(8)*y[8]+(16)*y[15]+(2)*y[16]+(12)*y[17]+(2)*y[18]+(12)*y[19]
+(22)*y[20]+(-14)*y[28]+(8)*y[38]+(14)*y[46];                              
(8)*v.128+(16)*v.135+(2)*v.136+(12)*v.137+(2)*v.138+(12)*v.139
+(22)*v.140+(-14)*v.148+(8)*v.158+(14)*v.166
> e*f;
(24)*v.241+(36)*v.242+(48)*v.243+(72)*v.244+(58)*v.245+(44)*v.246
+(30)*v.247+(16)*v.248
>h:=e*f;;
\end{verbatim}
Il faut pr{\'e}ciser pour cette orbite qu'avec les conventions de
\cite{Triplets}, l'{\'e}l{\'e}ment $X_{47}$ correspond {\`a} l'{\'e}l{\'e}ment {\tt
  x[46]=v.46} du logiciel, c'est pourquoi la d{\'e}finition du
$\mathfrak{sl}_{2}$-triplet est bien en accord avec \cite{Triplets}.

Calcul de {\tt g}, {\tt Bg}, {\tt z} et {\tt Bz}\,:
\begin{verbatim}
> g:=LieCentralizer(L,Subspace(L,[e]));
<Lie algebra of dimension 28 over Rationals>
> Bg:=BasisVectors(Basis(g));;
> z:=LieCentre(g);
<two-sided ideal in <Lie algebra of dimension 28 over Rationals>, 
  (dimension 3)>
> Bz:=BasisVectors(Basis(z));
[ v.15+v.16+v.17+v.18+v.19+v.20+v.38+v.46, v.119, v.120 ]
> h*Bz[2];
(14)*v.119
> h*Bz[3];
(16)*v.120
\end{verbatim}  
Les poids de {\tt z} sont 2,14,16; d'o{\`u} $m_{r}=16$. Il n'y a qu'une matrice {\`a} {\'e}tudier.

\begin{enumerate}
\item $m_{i_2}=14$,  $m_{k(2)}=4$
\begin{verbatim}
> h*Bg[3];
(4)*v.37
> h*Bg[6];
(4)*v.42+(-8)*v.44+(4)*v.45+(4)*v.48+(-4)*v.53+(4)*v.55
> h*Bg[7];
(4)*v.32+(2)*v.47+(2)*v.51+(2)*v.52+(4)*v.59
> h*Bg[8];
(4)*v.36+(-4)*v.40+(-8)*v.49+(-4)*v.54+(4)*v.57+(-4)*v.64+(4)*v.66

> ((f*Bg[8])*Bz[2]);
(-28)*v.120
\end{verbatim}
\end{enumerate}
{\bf Conclusion\,:} Cette orbite v{\'e}rifie la propri{\'e}t{\'e} $(P)$.\\

\item Caract{\'e}ristique\,:
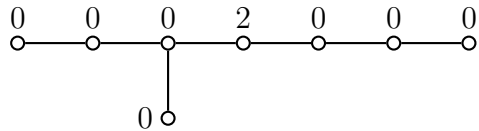
\begin{center}
\begin{pspicture}(-5,-1.2)(10,1)

\pscircle(-2,0){1mm}
\pscircle(-1,0){1mm}
\pscircle(0,0){1mm}
\pscircle(1,0){1mm}
\pscircle(2,0){1mm}
\pscircle(3,0){1mm}
\pscircle(4,0){1mm}

\pscircle(0,-1){1mm}

\psline(-1.9,0)(-1.1,0)
\psline(-0.1,0)(-0.9,0)
\psline(0.1,0)(0.9,0)
\psline(1.1,0)(1.9,0)
\psline(2.1,0)(2.9,0)
\psline(3.1,0)(3.9,0)
\psline(0,-0.1)(0,-0.9)

\rput[b](-2,0.2){$0$}
\rput[b](-1,0.2){$0$}
\rput[b](0,0.2){$0$}
\rput[b](1,0.2){$2$}
\rput[b](2,0.2){$0$}
\rput[b](3,0.2){$0$}
\rput[b](4,0.2){$0$}
\rput[r](-0.2,-1){$0$}
\end{pspicture}
\end{center}
D{\'e}finition du $\mathfrak{sl}_{2}$-triplet\,:
\begin{verbatim}
> e:=x[12]+x[21]+x[30]+x[31]+x[33]+x[42]+x[43]+x[53];                      
v.12+v.21+v.30+v.31+v.33+v.42+v.43+v.53
> f:=(5)*y[12]+y[21]+(5)*y[30]+(2)*y[31]+(8)*y[33]+(2)*y[42]+(8)*y[43]
+(9)*y[53];
(5)*v.132+v.141+(5)*v.150+(2)*v.151+(8)*v.153+(2)*v.162+(8)*v.163+(9)*v.173
> e*f;
(16)*v.241+(24)*v.242+(32)*v.243+(48)*v.244+(40)*v.245+(30)*v.246
+(20)*v.247+(10)*v.248
> h:=e*f;;
\end{verbatim}
Ici encore, il faut pr{\'e}ciser que l'{\'e}l{\'e}ment $X_{32}$ correspond 
{\`a} l'{\'e}l{\'e}ment {\tt x[31]=v.31} du logiciel et la d{\'e}finition du
$\mathfrak{sl}_{2}$-triplet est bien en accord avec \cite{Triplets}.

Calcul de {\tt g}, {\tt Bg}, {\tt z} et {\tt Bz}\,:
\begin{verbatim}
> g:=LieCentralizer(L,Subspace(L,[e]));
<Lie algebra of dimension 40 over Rationals>
> Bg:=BasisVectors(Basis(g));;
> z:=LieCentre(g);
<two-sided ideal in <Lie algebra of dimension 40 over Rationals>, 
  (dimension 5)>
> Bz:=BasisVectors(Basis(z));;
\end{verbatim}  
Les poids de {\tt z} sont 2,10,10,10,10; d'o{\`u} $m_{r}=10$. Il n'y a qu'une matrice {\`a}
{\'e}tudier.

\begin{enumerate}
\item $m_{i_2}=10$,  $m_{k(2)}=2$
\begin{verbatim}
> ((f*Bg[1])*Bz[2]); ((f*Bg[2])*Bz[2]); ((f*Bg[3])*Bz[2]); 
((f*Bg[4])*Bz[2]); ((f*Bg[5])*Bz[2]); ((f*Bg[6])*Bz[2]); 
((f*Bg[7])*Bz[2]); ((f*Bg[8])*Bz[2]); 
((f*Bg[9])*Bz[2]); ((f*Bg[10])*Bz[2]);
0*v.1
0*v.1
0*v.1
(-5)*v.120
(-10)*v.117
0*v.1
0*v.1
(-5)*v.118
0*v.1
(-10/3)*v.119
> ((f*Bg[1])*Bz[3]); ((f*Bg[2])*Bz[3]); ((f*Bg[3])*Bz[3]); 
((f*Bg[4])*Bz[3]); ((f*Bg[5])*Bz[3]); ((f*Bg[6])*Bz[3]); 
((f*Bg[7])*Bz[3]); ((f*Bg[8])*Bz[3]); 
((f*Bg[9])*Bz[3]); ((f*Bg[10])*Bz[3]);                     
0*v.1
(-2)*v.118
(-2)*v.120
(-2)*v.117
(-8)*v.118
0*v.1
0*v.1
v.120
(4/3)*v.119
0*v.1
> ((f*Bg[1])*Bz[4]); ((f*Bg[2])*Bz[4]); ((f*Bg[3])*Bz[4]); 
((f*Bg[4])*Bz[4]); ((f*Bg[5])*Bz[4]); ((f*Bg[6])*Bz[4]); 
((f*Bg[7])*Bz[4]); ((f*Bg[8])*Bz[4]); 
((f*Bg[9])*Bz[4]); ((f*Bg[10])*Bz[4]);
(-1)*v.119
0*v.1
0*v.1
0*v.1
(-9)*v.119
0*v.1
(-2)*v.118
0*v.1
(-2)*v.120
(2)*v.117
> ((f*Bg[1])*Bz[5]); ((f*Bg[2])*Bz[5]); ((f*Bg[3])*Bz[5); 
((f*Bg[4])*Bz[5]); ((f*Bg[5])*Bz[5]); ((f*Bg[6])*Bz[5]); 
((f*Bg[7])*Bz[5]); ((f*Bg[8])*Bz[5]); 
((f*Bg[9])*Bz[5]); ((f*Bg[10])*Bz[5]);
0*v.1
(-2)*v.120
0*v.1
v.118
(-8)*v.120
(-2)*v.118
(4/3)*v.119
(-2)*v.117
0*v.1
0*v.1
\end{verbatim}
La matrice {\`a} {\'e}tudier est de taille $4 \times 10$\,:
$$\left[
\begin{array}{cccccccccc}
0 & 0 & 0 & -2 \beta &-10\alpha &0 &0 &-2 \delta &0 & 2 \gamma \\
0 & -2\beta & 0 & \delta & -8 \beta & -2 \delta & -2 \gamma  & -5 \alpha & 0
& 0\\
-\gamma &0 & 0& 0&-9\gamma &0 &4/3 \delta & 0& 4/3 \beta & -10/3 \alpha \\
0 &-2 \delta & -2\beta&-5 \alpha &-8 \delta & 0&0 & \beta &-2 \gamma & 0
\end{array}
\right] \cdot
$$
Une {\'e}tude {\'e}l{\'e}mentaire permet de voir que cette matrice est de
rang 4, pour tout 4-uplet $(\alpha,\beta,\gamma,\delta)$ non nul.
\end{enumerate} 
{\bf Conclusion\,:} Cette orbite v{\'e}rifie la propri{\'e}t{\'e} $(P)$.\\
\end{enumerate}
{\bf Conclusion pour $E_8$\,:} Toutes les orbites nilpotentes distingu{\'e}es non
r{\'e}guli{\`e}res de $E_{8}$ v{\'e}rifient la propri{\'e}t{\'e} $(P)$.
 
\subsection{Calculs pour $F_{4}$}
\noindent D{\'e}finition de {\tt L}\,:
\begin{verbatim}
> L:=SimpleLieAlgebra("F",4,Rationals);
<Lie algebra of dimension 52 over Rationals>
> R:=RootSystem(L);
<root system of rank 4>
> P:=PositiveRoots(R);; 
> x:=PositiveRootVectors(R);
[ v.1, v.2, v.3, v.4, v.5, v.6, 
v.7, v.8, v.9, v.10, v.11, v.12, 
v.13, v.14, v.15, v.16, v.17, v.18, 
v.19, v.20, v.21, v.22, v.23, v.24 ]
> y:=NegativeRootVectors(R);
[ v.25, v.26, v.27, v.28, v.29, v.30, 
v.31, v.32, v.33, v.34, v.35, v.36, 
v.37, v.38, v.39, v.40, v.41, v.42, 
v.43, v.44, v.45, v.46, v.47, v.48 ]
> CanonicalGenerators(R)[3];
[ v.49, v.50, v.51, v.52 ]
\end{verbatim}
Dans $F_{4}$, il y a trois orbites nilpotentes distingu{\'e}es non
r{\'e}guli{\`e}res. Pour $F_{4}$, les conventions du logiciel GAP4 sont tr{\`e}s
diff{\'e}rentes de celles adopt{\'e}es dans \cite{Triplets}; dans \cite{Triplets}, le
diagramme de Dynkin est\,:
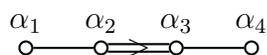
\begin{center}
\begin{pspicture}(-5,-0.2)(10,1)

\pscircle(-2,0){1mm}
\pscircle(-1,0){1mm}
\pscircle(0,0){1mm}
\pscircle(1,0){1mm}

\psline(-1.9,0)(-1.1,0)
\psline(-0.09,0.05)(-0.91,0.05)
\psline(-0.09,-0.05)(-0.91,-0.05)
\psline(0.1,0)(0.9,0)

\rput[b](-2,0.2){$\alpha_{1}$}
\rput[b](-1,0.2){$\alpha_{2}$}
\rput[b](0,0.2){$\alpha_{3}$}
\rput[b](1,0.2){$\alpha_{4}$}

\rput(-0.5,0){$>$}
\end{pspicture}
\end{center}
Il semble que dans GAP4 le diagramme de Dynkin soit plut{\^o}t le suivant\,:
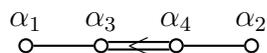
\begin{center}
\begin{pspicture}(-5,-0.2)(10,1)

\pscircle(-2,0){1mm}
\pscircle(-1,0){1mm}
\pscircle(0,0){1mm}
\pscircle(1,0){1mm}

\psline(-1.9,0)(-1.1,0)
\psline(-0.09,0.05)(-0.91,0.05)
\psline(-0.09,-0.05)(-0.91,-0.05)
\psline(0.1,0)(0.9,0)

\rput[b](-2,0.2){$\alpha_{1}$}
\rput[b](-1,0.2){$\alpha_{3}$}
\rput[b](0,0.2){$\alpha_{4}$}
\rput[b](1,0.2){$\alpha_{2}$}

\rput(-0.5,0){$<$}
\end{pspicture}
\end{center}
Par cons{\'e}quent, il est difficile d'utiliser directement les donn{\'e}es de
\cite{Triplets} dans GAP4. On utilise les correspondances suivantes\,:
{\tt x[1]}=$X_{4}$, {\tt x[2]}=$X_{1}$, {\tt x[3]}=$X_{3}$, {\tt x[4]}=$X_{1}$, {\tt x[5]}=$X_{7}$, {\tt
  x[6]}=$X_{5}$, {\tt x[6]}=$X_{5}$, {\tt x[7]}=$X_{6}$, {\tt
  x[8]}=$X_{10}$, {\tt x[9]}=$X_{8}$,  {\tt x[10]}=$X_{9}$,  {\tt
  x[18]}=$X_{18}$. Cependant, m{\^e}me avec ces relations, les
$\mathfrak{sl}_{2}$-triplets de \cite{Triplets} ne conviennent pas. On
utilise la commande {\tt FindSl2} qui permet de chercher une
sous-alg{\`e}bre {\tt s} isomorphe {\`a} $\mathfrak{sl}_{2}$ et contenant {\tt
  e}. On s'assure auparavant que l'{\'e}l{\'e}m{\'e}nt {\tt e} est bien nilpotent {\`a}
l'aide de la commande {\tt IsNilpotentElement} et on v{\'e}rifie aussi que
le $\mathfrak{sl}_{2}$-triplet obtenu correspond bien {\`a} la
caract{\'e}ristique voulue.\\

\begin{enumerate}
\item Caract{\'e}ristique\,:
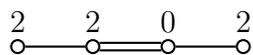
\begin{center}
\begin{pspicture}(-5,-0.2)(10,1)

\pscircle(-2,0){1mm}
\pscircle(-1,0){1mm}
\pscircle(0,0){1mm}
\pscircle(1,0){1mm}

\psline(-1.9,0)(-1.1,0)
\psline(-0.09,0.05)(-0.91,0.05)
\psline(-0.09,-0.05)(-0.91,-0.05)
\psline(0.1,0)(0.9,0)

\rput[b](-2,0.2){$2$}
\rput[b](-1,0.2){$2$}
\rput[b](0,0.2){$0$}
\rput[b](1,0.2){$2$}
\end{pspicture}
\end{center}
D{\'e}finition du $\mathfrak{sl}_{2}$-triplet\,:
\begin{verbatim}
> e:=x[2]+x[4]+x[5]+x[7];
v.2+v.4+v.5+v.7
> IsNilpotentElement(L,e);                   
true
> s:=FindSl2(L,a);                           
<Lie algebra of dimension 3 over Rationals>
> Bs:=BasisVectors(Basis(s));                
[ v.2+v.4+v.5+v.7, v.49+(7/5)*v.50+(9/5)*v.51+(13/5)*v.52, 
  v.25+(7/5)*v.26+v.28+v.29+(4/5)*v.31+(-4/5)*v.34 ]
\end{verbatim}
L'{\'e}l{\'e}ment central de cette base est {\'e}gal au dixi{\`e}me de l'{\'e}lement neutre correspondant {\`a} la caract{\'e}ristique; par suite
en prenant pour {\tt f} dix fois le troisi{\`e}me {\'e}l{\'e}ment de cette
base, on obtient un  $\mathfrak{sl}_{2}$-triplet pour cette caract{\'e}ristique.
\begin{verbatim}
> h:=(14)*H[2]+(26)*H[4]+(18)*H[3]+(10)*H[1];
(10)*v.49+(14)*v.50+(18)*v.51+(26)*v.52
> f:=(10)*y[1]+(14)*y[2]+(10)*y[4]+(10)*y[5]+(8)*y[7]+(-8)*y[10];
(10)*v.25+(14)*v.26+(10)*v.28+(10)*v.29+(8)*v.31+(-8)*v.34
> e*f;
(10)*v.49+(14)*v.50+(18)*v.51+(26)*v.52
\end{verbatim}
Calcul de {\tt g}, {\tt Bg}, {\tt z} et {\tt Bz}\,:
\begin{verbatim}
> g:=LieCentralizer(L,Subspace(L,[e]));
<Lie algebra of dimension 6 over Rationals>
> Bg:=BasisVectors(Basis(g));;
> z:=LieCentre(g);
<two-sided ideal in <Lie algebra of dimension 6 over Rationals>, 
(dimension 3 )>
> Bz:=BasisVectors(Basis(z));
[ v.2+v.4+v.5+v.7, v.20+v.21+(2)*v.22, v.24 ]
> h*Bz[2];
(10)*v.20+(10)*v.21+(20)*v.22
> h*Bz[3];
(14)*v.2
\end{verbatim}  
Les poids de {\tt z} sont 2,10,14; d'o{\`u} $m_{r}=14$. Il n'y a qu'une matrice {\`a} {\'e}tudier.

\begin{enumerate}
\item $m_{i_2}=10$,  $m_{k(2)}=6$
\begin{verbatim}
> h*Bg[3];
(6)*v.11+(6)*v.14+(-12)*v.15+(6)*v.16

> ((f*Bg[3])*Bz[2]);
(30)*v.24
\end{verbatim}
\end{enumerate}
{\bf Conclusion\,:} Cette orbite v{\'e}rifie la propri{\'e}t{\'e} $(P)$.\\

\item Caract{\'e}ristique\,:
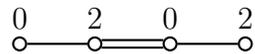
\begin{center}
\begin{pspicture}(-5,-0.2)(10,1)

\pscircle(-2,0){1mm}
\pscircle(-1,0){1mm}
\pscircle(0,0){1mm}
\pscircle(1,0){1mm}

\psline(-1.9,0)(-1.1,0)
\psline(-0.09,0.05)(-0.91,0.05)
\psline(-0.09,-0.05)(-0.91,-0.05)
\psline(0.1,0)(0.9,0)

\rput[b](-2,0.2){$0$}
\rput[b](-1,0.2){$2$}
\rput[b](0,0.2){$0$}
\rput[b](1,0.2){$2$}
\end{pspicture}
\end{center}
D{\'e}finition du $\mathfrak{sl}_{2}$-triplet\,:
\begin{verbatim}
> e:=x[1]+x[6]+x[5]+x[10];
v.1+v.5+v.6+v.10
> IsNilpotentElement(L,e);
true
> s:=FindSl2(L,e);
<Lie algebra of dimension 3 over Rationals>
> Bs:=BasisVectors(Basis(s));
[ v.1+v.5+v.6+v.10, v.49+(5/4)*v.50+(7/4)*v.51+(5/2)*v.52, 
  v.25+v.29+(5/2)*v.30+v.31+v.33+(5/2)*v.34 ]
\end{verbatim}
L'{\'e}l{\'e}ment central de cette base est {\'e}gal au huiti{\`e}me de l'{\'e}l{\'e}ment neutre de
la carat{\'e}ristique. Ici {\tt Bs[1]*Bs[3]=2*Bs[2]}; on prend alors pour
{\tt f} quatre fois le troisi{\`e}me {\'e}l{\'e}ment de cette base. 
\begin{verbatim}
> h:=(10)*H[2]+(20)*H[4]+(14)*H[3]+(8)*H[1];
(8)*v.49+(10)*v.50+(14)*v.51+(20)*v.52
> f:=(4)*y[1]+(4)*y[5]+(10)*y[6]+(4)*y[7]+(4)*y[9]+(10)*y[10];
(4)*v.25+(4)*v.29+(10)*v.30+(4)*v.31+(4)*v.33+(10)*v.34
> e*f;
(8)*v.49+(10)*v.50+(14)*v.51+(20)*v.52
\end{verbatim}
Calcul de {\tt g}, {\tt Bg}, {\tt z} et {\tt Bz}\,:
\begin{verbatim}
> g:=LieCentralizer(L,Subspace(L,[e]));
<Lie algebra of dimension 8 over Rationals>
> Bg:=BasisVectors(Basis(g));;
> z:=LieCentre(g);
<two-sided ideal in <Lie algebra of dimension 8 over Rationals>, 
(dimension 3)>
> Bz:=BasisVectors(Basis(z));
[ v.1+v.5+v.6+v.10, v.23, v.24 ]
> h*Bz[2];
(10)*v.23
> h*Bz[3];
(10)*v.24
\end{verbatim}  
Les poids de {\tt z} sont 2,10,10; d'o{\`u} $m_{r}=10$. Il n'y a qu'une matrice {\`a} {\'e}tudier.

\begin{enumerate}
\item $m_{i_2}=10$,  $m_{k(2)}=2$
\begin{verbatim}
> h*Bg[2];
(2)*v.4+(2)*v.5+(2)*v.6+(-2)*v.7+(2)*v.10
> h*Bg[3];
(2)*v.5+(2)*v.9+(-2)*v.13
> h*Bg[4];
(4)*v.8+(4)*v.11+(-4)*v.12+(-4)*v.14+(-4)*v.16

> ((f*Bg[1])*Bz[2]); ((f*Bg[2])*Bz[2]); ((f*Bg[3])*Bz[2]);
(-8)*v.23
(-2)*v.23
(-2)*v.24
> ((f*Bg[1])*Bz[3]); ((f*Bg[2])*Bz[3]); ((f*Bg[3])*Bz[3]);
(-2)*v.23
(2)*v.23+(-10)*v.24
(-8)*v.24
\end{verbatim}
La matrice {\`a} {\'e}tudier est 
$$\left[
\begin{array}{ccc}
-8 \alpha -2 \beta & -2 \alpha +2 \beta & 0          \\
0 & -10 \beta  & -2 \alpha -8 \beta 
\end{array}
\right] \cdot $$ 
Une br{\`e}ve {\'e}tude de cette matrice montre qu'elle est de rang 2 pour
tout couple $(\alpha,\beta)$ non
nul. 
\end{enumerate}
{\bf Conclusion\,:} Cette orbite v{\'e}rifie la propri{\'e}t{\'e} $(P)$.\\

\item Caract{\'e}ristique\,:
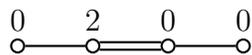
\begin{center}
\begin{pspicture}(-5,-0.2)(10,1)

\pscircle(-2,0){1mm}
\pscircle(-1,0){1mm}
\pscircle(0,0){1mm}
\pscircle(1,0){1mm}

\psline(-1.9,0)(-1.1,0)
\psline(-0.09,0.05)(-0.91,0.05)
\psline(-0.09,-0.05)(-0.91,-0.05)
\psline(0.1,0)(0.9,0)

\rput[b](-2,0.2){$0$}
\rput[b](-1,0.2){$2$}
\rput[b](0,0.2){$0$}
\rput[b](1,0.2){$0$}
\end{pspicture}
\end{center}
D{\'e}finition du $\mathfrak{sl}_{2}$-triplet\,:
\begin{verbatim}
> e:=x[9]+x[10]+x[8]+x[18];   
v.8+v.9+v.10+v.18
> f:=2*y[9]+2*y[10]+2*y[8]+2*y[18];
(2)*v.32+(2)*v.33+(2)*v.34+(2)*v.42
> e*f;
(4)*v.49+(6)*v.50+(8)*v.51+(12)*v.52
> h:=e*f;;
\end{verbatim}
Pour cette orbite, les donn{\'e}es de \cite{Triplets} conviennent, pour
des raisons qui m'{\'e}chappent.\\
\\
Calcul de {\tt g}, {\tt Bg}, {\tt z} et {\tt Bz}\,:
\begin{verbatim}
> g:=LieCentralizer(L,Subspace(L,[e]));
<Lie algebra of dimension 12 over Rationals>
> Bg:=BasisVectors(Basis(g));;
> z:=LieCentre(g);
<two-sided ideal in <Lie algebra of dimension 12 over Rationals>, 
(dimension 3 )>
> Bz:=BasisVectors(Basis(z));
[ v.8+v.9+v.10+v.18, v.23, v.24 ]
> h*Bz[2];
(6)*v.23
> h*Bz[3];
(6)*v.24
\end{verbatim}  
Les poids de {\tt z} sont 2,6,6; d'o{\`u} $m_{r}=6$. Il n'y a qu'une matrice {\`a} {\'e}tudier.

\begin{enumerate}
\item $m_{i_2}=16$,  $m_{k(2)}=2$
\begin{verbatim}
> h*Bg[1]; h*Bg[2]; h*Bg[3]; h*Bg[4]; h*Bg[5];
(2)*v.8+(2)*v.9
(2)*v.4+(-1)*v.11+v.12
(2)*v.12+(-2)*v.13
(2)*v.6+(-1)*v.7+v.14
(2)*v.14+(-2)*v.15
> h*Bg[7];                                                                 
(2)*v.10+(2)*v.18

> ((f*Bg[1])*Bz[2]); ((f*Bg[2])*Bz[2]); ((f*Bg[3])*Bz[2]);
((f*Bg[4])*Bz[2]); ((f*Bg[5])*Bz[2]); ((f*Bg[7])*Bz[2])
(-4)*v.23
(-2)*v.24
(-2)*v.24
0*v.1
0*v.1
(-2)*v.23

> ((f*Bg[1])*Bz[3]); ((f*Bg[2])*Bz[3]); ((f*Bg[3])*Bz[3]);
((f*Bg[4])*Bz[3]); ((f*Bg[5])*Bz[3]); ((f*Bg[7])*Bz[3])
(-4)*v.24
0*v.1
0*v.1
(-2)*v.23
(-2)*v.23
(-2)*v.24
\end{verbatim}
La matrice {\`a} {\'e}tudier est 
$$\left[
\begin{array}{cccccc}
-4 \alpha & 0          & 0           & -2 \beta  & -2 \beta & -2 \alpha \\
-4 \beta  & -2 \alpha  & -2 \alpha  & 0 & 0 & -2 \beta
\end{array}
\right] \cdot $$ 
C'est une matrice de rang 2 pour tout couple $(\alpha,\beta)$ non
nul. 
\end{enumerate}
{\bf Conclusion\,:} Cette orbite v{\'e}rifie la propri{\'e}t{\'e} $(P)$.\\
\end{enumerate}
{\bf Conclusion pour $F_{4}$\,:} Toutes les orbites nilpotentes distingu{\'e}es non
r{\'e}guli{\`e}res de $F_{4}$ v{\'e}rifient la propri{\'e}t{\'e} $(P)$.

\nocite{*}
\bibliographystyle{plain}
\bibliography{annexe}

\vspace{3cm}

\normalsize
~\\
{\sc Anne Moreau \\
~\\
Universit{\'e} Paris 7 - Denis Diderot,\\
Institut de Math{\'e}matiques de Jussieu,\\
Th{\'e}orie des groupes,\\
Case 7012\\
2, Place jussieu\\
75251 Paris Cedex 05, France. }\\
~\\
{\em E-mail\,: {\tt moreaua@math.jussieu.fr} }

\end{document}

%% file: francais_bis.tex
\usepackage{pstricks,pst-plot}
\usepackage{amsfonts}

\usepackage[latin1]{inputenc}
\usepackage[frenchb]{babel}  
\def\g{\mathfrak{g}}

\def\z{\mathfrak{z}}

\def\n{\mathfrak{n}}

\def\sl2{\mathfrak{sl}_{2}-{\rm triplet}}

\let \epsilon=\varepsilon
\def \phi {\varphi}

